\documentclass[reqno, 12pt]{amsart}
\usepackage{amsfonts}
\usepackage{amssymb,amsmath}
\usepackage{amsthm}
\usepackage{amscd}
\usepackage{upref}

\makeatletter
\def\LaTeX{\leavevmode L\raise.42ex
    \hbox{\kern-.3em\size{\sf@size}{0pt}\selectfont A}\kern-.15em\TeX}
 \makeatother

\numberwithin{equation}{section}

\newtheorem{lemma}{Lemma}[section]
\newtheorem{theorem}[lemma]{Theorem} 
\newtheorem{corollary}[lemma]{Corollary}
\newtheorem{proposition}[lemma]{Proposition}

\theoremstyle{definition}

\newtheorem{definition}[lemma]{Definition}

\renewcommand{\det}{\operatorname{Det}}

  \newcommand{\e}{\eqref}

\newcommand{\q}{\quad}

\newcommand{\ti}{\tilde}
\newcommand{\wt}{\widetilde}

\newcommand{\la}{\langle}
\newcommand{\ra}{\rangle}

\newcommand{\ov}{\overline}
 \renewcommand{\d}{\delta}

   \newcommand{\sgn}{\operatorname{sgn}}

  \newcommand{\rank}{\operatorname{rank}}

\renewcommand\Im{\operatorname{Im}}
\renewcommand\Re{\operatorname{Re}}

\newenvironment{pf}{\begin{proof}}{\end{proof}}

\def\qqq{\mathrel{\subset\mkern-15mu\lower.38ex\hbox{${\scriptscriptstyle\rightarrow}$}}}
\let\cal\mathcal

\let\Bbb\mathbb

    \DeclareMathOperator{\diag}{diag}
  
\begin{document}

\title{On finite rank Hankel operators}
\author{ D. R. Yafaev}
\address{ IRMAR, Universit\'{e} de Rennes I\\ Campus de
  Beaulieu, 35042 Rennes Cedex, FRANCE}
\email{yafaev@univ-rennes1.fr}
\keywords{Hankel  operators,  operators of finite rank, the sign-function,   necessary and sufficient conditions for    the  sign-definiteness,    total multiplicity of the positive and negative spectra,  the Carleman operator and its perturbations}
\subjclass[2000]{47A40, 47B25}

%  \date{\today}

\begin{abstract}
  For self-adjoint Hankel  operators of   finite rank, we  find an explicit formula for the total multiplicity of their negative and positive spectra. We also show that very strong perturbations, for example,   a  perturbation  by the Carleman operator, do not change the total  number of negative eigenvalues of finite rank  Hankel operators.
  \end{abstract}

\maketitle

% \thispagestyle{empty}

%************************************************************
\section{Introduction. Main results}  
%***********************************************************

{\bf 1.1.}
Hankel operators can be defined as integral operators
\begin{equation}
(H f)(t) = \int_{0}^\infty h(t+s) f(s)ds 
\label{eq:H1}\end{equation}
in the space $L^2 ({\Bbb R}_{+}) $ with kernels $h$ that depend  on the sum of variables only. Of course $H$ is symmetric if $  h(t)=\ov{h(t)}$.  

Integral kernels of self-adjoint Hankel operators $H$ of  finite rank are given (this is the Kronecker theorem -- see, e.g., Sections~1.3 and 1.8 of the book \cite{Pe}) by the formula
\begin{equation}
h(t) =\sum_{m=1}^M  P_{m} (t)e^{- \alpha_{m}  t}
\label{eq:FDvm}\end{equation}
where $\Re\alpha_{m}  >0$ and 
 $   P_{m} (t)$ 
    are polynomials of degree $K_{m} $. If $H$ is self-adjoint,  then necessarily the sum in  \e{eq:FDvm} contains both exponentials  $e^{- \alpha_{m}  t}$ and   $e^{- \bar{\alpha}_{m}  t}$. Let $\Im \alpha_{m}=0$ for $m=1,\ldots, M_{0}$, $\Im \alpha_{m}>0$   and $\alpha_{ M_{1}+m}=\bar{\alpha}_{ m}$ for $m=M_{0} +1, \ldots,M_{0} + M_{1} $. Thus $M=M_{0}+2M_{1}$; of course the cases $M_{0}=0$ or $M_{1}=0$ are not excluded. The condition  $ h(t)=\ov{h(t)} $ requires also that $ P_{m} (t) = \ov{P_{m} (t)} $ for $m=1,\ldots, M_{0}$ and $ P_{M_{1}+m} (t) = \ov{P_{m} (t)}$  for $m=M_{0} +1, \ldots,M_{0} + M_{1} $.    As is well known and as we shall see below,
      \[
\rank H=\sum_{m=1}^M K_{m}+M=: r.
\]

     For $m=1,\ldots, M_{0}$, we
 set
\begin{equation}
 {\sf p}_{m}=P_{m}^{(K_m)},
\label{eq:ppq}\end{equation}
that is, ${\sf p}_{m}/K_{m}!$ is the coefficient at $t^{K_{m}}$ in the polynomial  $P_{m} (t)$, 
 and
    \begin{equation}
\left. \begin{aligned}
{\cal N}_{+}^{(m)}  = {\cal N}_{-}^{(m)} &=(K_{m} +1)/2 \q {\rm if} \; K_{m} \; {\rm is  } \; {\rm   odd}
 \\
{\cal N}_{+}^{(m)}  -1= {\cal N}_{-}^{(m)}  &= K_{m}  /2 \q {\rm if} \; K_{m} \; {\rm is  } \; {\rm   even} \; {\rm   and}
 \q {\sf p}_{m}  > 0
 \\
{\cal N}_{+}^{(m)}= {\cal N}_{-}^{(m)}  -1  &= K_{m}  /2  \q {\rm if} \; K_{m} \; {\rm is  } \; {\rm   even} \; {\rm   and}
 \q {\sf p}_{m} < 0 .
    \end{aligned}
    \right\}
    \label{eq:RXm}\end{equation}

 For a self-adjoint operator $A$, we denote by $N_{+}(A)$ (by $N_{-}(A)$) the total mutiplicity of its strictly positive (negative) spectrum.
   Our main result is formulated as follows.
 
  \begin{theorem}\label{FDH1}
 Let $H$ be the self-adjoint Hankel operator  of  finite rank with kernel $h(t)$ given by formula \e{eq:FDvm} where $P_{m} (t)$ are polynomials of degree $K_{m}$, and let the numbers ${\cal N}_\pm^{(m)}$ be defined by formula \e{eq:RXm}.  
  Then the total numbers $N_{\pm} (H)$  of   $($strictly$)$ positive and negative eigenvalues of the operator $ H$ are given by the formula 
   \begin{equation}
 N_{\pm} (H) =  \sum_{m=1}^{M_{0}}  {\cal N}_\pm^{(m)} + \sum_{m=M_{0}+1}^{M_{0}+M_1}  K_{m}  +M_{1} .
   \label{eq:TN}\end{equation}
 \end{theorem}

 Formula \e{eq:TN} shows that every pair of complex conjugate terms 
     \begin{equation}
     P_{m} (t) e^{-\alpha_{m}t}+ \ov{P_{m } (t)} e^{-\bar{\alpha}_{m }t}, \q m= M_{0}+1,\ldots, M_{0}+M_{1},
      \label{eq:TNX}\end{equation}
  in representation \e{eq:FDvm} of $h(t)$  
 yields $K_{m}+1$ positive and $K_{m}+1$ negative eigenvalues. In view of \e{eq:RXm} the contribution of every real term
 $P_{m} (t)e^{- \alpha_{m}  t}$ also consists of the equal numbers $(K_{m}+1)/2$ of positive and  negative eigenvalues if the degree $K_{m}$ of $P_m (t)$ is odd. If $K_{m}$   is even, then there is one extra positive (negative) eigenvalue if $P_{m}^{(K_m)}>0$
 ($P_{m}^{(K_m)}<0$).  In particular,  in the question considered, there is no ``interference" between different real terms $P_{m} (t) e^{-\alpha_{m}t}$, $m=1,\ldots, M_{0}$, and pairs \e{eq:TNX}.
   
According to \e{eq:TN} the operator $H$ cannot be sign-definite if $M_{1}>0$.  Moreover, 
according to \e{eq:RXm} the operator $H$ cannot be sign-definite if $K_m >0$ at least for one $m=1,\ldots, M_{0}$.
Therefore we have the following result.
    
  \begin{corollary}\label{FD}
A  Hankel operator  $H$ of finite rank in the space $L^2 ({\Bbb R}_{+})$ is positive  $($negative$)$ if
and only if its kernel is given by the formula 
   \[
h(t) =\sum_{m= 1}^{M_{0}} {\sf p}_{m} e^{- \alpha_{m} t}
\]
%\label{eq:FD1}\end{equation}
where $\alpha_{m}>0$ and $  {\sf p}_{m} >0$ $(  {\sf p}_{m} < 0)$.
 \end{corollary}
 
 Let us recall the paper \cite{MPT} by
A.~V.~Megretskii, V.~V.~Peller, and S.~R.~Treil. In the particular case of finite rank Hankel operators $H$, it follows from the results of \cite{MPT} that the spectra of $H$ are characterized by the condition that the multiplicities of eigenvalues $\lambda$ and $-\lambda$ do not differ by more than $1$. Compared to Theorem~\ref{FDH1}, this result is of a completely different nature.

 \medskip
 
 {\bf 1.2.}
  The result of Theorem~\ref{FDH1} turns out to be   stable under a large class of perturbations of finite rank Hankel operators. As an example, we consider the sum $H=H_{0}+V$
  of the Carleman operator $H_{0}$, that is, of
the Hankel operator with   kernel $h_{0}(t)=  t^{-1}$,  and of a finite rank Hankel operator $V$.  Recall that the Carleman operator  has the absolutely continuous spectrum $[0,\pi]$ of multiplicity $2$. We obtain the following result.

  \begin{theorem}\label{FDHC}
 Let     $ H_{0}$ be  the Hankel operator  with kernel $h_{0}(t)=t^{-1}$. If $V$ is a Hankel operator  of   finite rank and $H=H_{0}+V$,  then 
  $$N_{-} (H)=N_{-} (V).$$
  In particular, $H \geq 0$ if and only if $V\geq 0$.
 \end{theorem}

  The inequality $N_{-} (H)\leq N_{-} (V)$ is of course obvious because $H_{0}\geq 0$. On the contrary, the opposite inequality $N_{-} (H)\geq N_{-} (V)$ looks surprising because the Carleman operator is ``much stronger" than $V$; it is not even compact. Nevertheless its adding  does not change
    the total number of   negative   eigenvalues.

   It is natural to compare (this point of view goes back to J.~S.~Howland \cite{Howland})   Hankel operators $H$ with ``perturbed" kernels $h (t)=t^{-1}+v(t)$    to Schr\"odinger operators  $D^2+ {\sf V}(x)$. The assumption that
   $v(t)$    decays sufficiently rapidly as $t\to \infty$ and is not too singular as $t\to 0$ corresponds to a sufficiently rapid decay 
 of a potential $  {\sf V}(x)$  as $|x|\to \infty$.   As shown in \cite{Y2},   the results     on the discrete spectrum of the operator $H$ lying {\it above} its essential spectrum $[0,\pi]$ are close in spirit to the results on the discrete (negative) spectrum of the Schr\"odinger operator  $D^2+ {\sf V}(x)$. On the contrary,  according to Theorem~\ref{FDHC}   the results on the negative spectrum of  Hankel operators are drastically different from those for the Schr\"odinger operators.

 \medskip
 
 {\bf 1.3.} 
 Our proofs of Theorems~\ref{FDH1} and \ref{FDHC} rely on the approach suggested in 
 \cite{Y}. It is shown in  \cite{Y} that a Hankel operator $H$ has the same numbers of negative and positive eigenvalues as an operator $S$ of multiplication by some function $s(x)$. In particular, $\pm H\geq 0$ if and only if $\pm S\geq 0$. Therefore we use the term ``sign-function" for $s(x)$. In specific examples 
  functions $s(x)$ may be of a quite different nature. For instance, for finite rank Hankel operators,  $s(x)$ is  a distribution which is an explicit combination of delta functions and their derivatives. This  allows us to calculate  the total  numbers of negative and positive eigenvalues of such operators and thus prove Theorem~\ref{FDH1}.

  As far as Theorem~\ref{FDHC} is concerned, we note that
    the sign-function of the Carleman operator  equals $1$. Its support is essentially disjoint from supports of the sign-functions of  finite rank Hankel operators $V$. Very loosely speaking, it means  that the operators $H_{0}$ and $V$ ``live in orthogonal subspaces", and hence the positive operator $H_{0}$ does not affect the negative spectrum of $H=H_{0}+V$.
    
    We note that since the sign-function is, in general, a distribution,  $S$ need not be defined as an operator. Therefore we work with quadratic forms which is both more general and more convenient.
 
 Roughly speaking, the approach of  \cite{Y} can be described as follows.
  Let $S$ be the {\it formal} operator of multiplication by the function $s(x)$. Then   the identity  
   \begin{equation}
 H=  \Xi^* S \Xi
\label{eq:MAIDs}\end{equation} 
holds with some invertible operator $\Xi$. It follows that
  \begin{equation}
N_{\pm}(H)= N_{\pm}(S). 
\label{eq:MAIDT}\end{equation} 

For finite rank Hankel operators $H$, the form $(Su,u)$ is determined by values of functions $u(x)$ and their  derivatives at some finite number of points. Therefore number
\e{eq:MAIDT} equals $N_{\pm}({\bf S})$ for some Hermitian matrix $\bf S$ (the sign-matrix of the operator $H$) constructed in terms of $s(x)$. It turns out that the matrix $\bf S$ has a very special structure which allows us to calculate the number
 $N_{\pm}({\bf S})$ explicitly.
 
   \medskip
 
 {\bf 1.4.}
 Let us briefly describe the structure of the paper. We collect necessary results of \cite{Y} 
   in Section~2. Proofs of Theorems~\ref{FDH1} and \ref{FDHC} are given in Section~3.
    Hankel operators can be standardly realized not only in $L^2({\Bbb R}_{+})$ but also in the  Hardy spaces  ${\Bbb H}_{+}^2({\Bbb R})$, ${\Bbb H}_{+}^2({\Bbb T})$ and in the space of sequences $l^2({\Bbb Z}_{+})$. The interrelations between different representations are discussed in the auxiliary Section~4. This information is used in Section~5 to     reformulate Theorems~\ref{FDH1} and \ref{FDHC} in  the spaces ${\Bbb H}_{+}^2({\Bbb R})$, ${\Bbb H}_{+}^2({\Bbb T})$ and   $l^2({\Bbb Z}_{+})$. Finally, in the Appendix we describe the group of  automorphisms of the set of Hankel operators in all these spaces as well as in the space $L^2({\Bbb R}_{+})$.

  Let us introduce some standard
 notation. We first recall that ${\Bbb T}$ is the unit circle in the complex plane and ${\Bbb Z}_{+} $ is the set of all nonnegative integers.
  We denote by $\Phi$, 
\[
(\Phi u) (\xi)=  (2\pi)^{-1/2} \int_{-\infty}^\infty u(x) e^{ -i x \xi} dx,
\]
  the Fourier transform.   The space $\cal Z= \cal Z ({\Bbb R})$ of test functions  is defined as the subset  
 of the Schwartz  space ${\cal S}={\cal S} ({\Bbb R}) $ which consists of functions $\varphi  $ admitting the analytic continuation to   entire functions in the   complex plane $\Bbb C$   and satisfying,   for all $z\in \Bbb C$,   bounds 
 \[
  | \varphi (z)| \leq C_{n}  (1+| z |)^{-n} e^{r |\Im z |}
  \]
  for some $r=r(\varphi)>0$  and all $n$. We recall that the Fourier transform 
  $\Phi : \cal Z \to C_{0}^\infty ({\Bbb R})$
  and
  $\Phi^*: C_{0}^\infty ({\Bbb R}) \to \cal Z$. 
  
   The dual classes of distributions (continuous antilinear functionals) are denoted ${\cal S}'$, $C_{0}^\infty ({\Bbb R})'$ and ${\cal Z}'$, respectively. 
  We use the notation   ${\pmb\la} \cdot, \cdot {\pmb\ra}$ and $\la \cdot, \cdot\ra$ for   the  
  duality symbols in $L^2 ({\Bbb R}_{+})$ and $L^2 ({\Bbb R})$, respectively. They are   linear in the first argument and antilinear in the second argument.

   The Dirac function is standardly denoted $\d(\cdot)$; $\d_{n,m}$ is the Kronecker symbol, i.e., $\d_{n,n}=1$  and $\d_{n,m}= 0$  if $n\neq m$. 
   The letter $C$ (sometimes with indices) denotes various positive constants whose precise values are inessential.

%************************************************************
\section{The sign-function }  
%*

Here we briefly discuss necessary results of \cite{Y} adapting them to the case of bounded Hankel operators.

\medskip

{\bf 2.1.}
Let us consider a Hankel operator $H$ defined by equality \e{eq:H1} in the space $L^2 ({\Bbb R}_{+})$.   Actually, it is more convenient to work with sesquilinear forms  instead of operators. Let us introduce the   Laplace convolution
\begin{equation} 
( \bar{f}_{1}\star f_{2})(t)=
\int_{0}^t    \overline{f_{1}(s)} f_{2}(t-s) ds  
 \label{eq:HH1}\end{equation} 
 of   functions $ \bar{f}_{1}$ and $ f_{2}$. Then
\begin{equation}
(Hf_{1} ,f_{2})=  {\pmb \la} h, \bar{f}_{1}\star f_{2} {\pmb \ra}=: h[f_{1},f_{2}]
 \label{eq:HH}\end{equation}
 where we write ${\pmb \la} \cdot, \cdot {\pmb \ra}$ instead of $( \cdot, \cdot )$ because $h$ may be a distribution.
 
We consider form \e{eq:HH} on elements ${f}_{1},  f_{2} \in \cal D$ where $\cal D$ is defined as follows. Put  
\[
(U f)(x) =e^{x/2} f(e^x).
\]
Then    $U: L^2 ({\Bbb R}_{+})\to L^2 ({\Bbb R} )$  is the unitary operator.
The set $\cal D$ consists of functions $f  (t)$  such that $  Uf  \in {\cal Z}$. Since
\[
f(t)= t^{-1/2}  (Uf) (\ln t)
\]
and ${\cal Z}\subset {\cal S}$, we see that   functions $f  \in \cal D$ and their derivatives satisfy the estimates
 \[
|f^{(m)}(t)|=  C_{n,m} t^{-1/2-m} (1+|\ln t|)^{-n}
\]
for all $n$ and $m$. Of course, $\cal D$  is dense in the space $L^2 ({\Bbb R}_{+})$.
It is shown in \cite{Y} that  if  ${f}_{1},  f_{2} \in \cal D$, then  the function
 \[
\Omega(x)=(\bar{f}_{1}\star f_{2} )(e^x) 
\]
 belongs to the set  ${\cal Z}$.

  With respect to $h$, we assume that the distribution  
 \begin{equation}
 \theta (x)= e^x h(e^x)
 \label{eq:HH3}\end{equation}
 is an element of the space $\cal Z'$. The set of all such  $h$ will be denoted   $ {\cal Z}'_{+}$,   that is,
 \[
h\in  {\cal Z}_{+}' \Longleftrightarrow \theta\in  {\cal Z}'.
\]
It is shown in \cite{Y} that {\it this condition is satisfied for all bounded Hankel operators }$H$. Since $\Omega\in {\cal Z}$, the form
 \[
  {\pmb\la} h, \bar{f}_{1}\star f_{2}  {\pmb\ra} = \int_{0}^\infty h( t ) ( {f}_{1}\star \bar{f}_{2} )(t) d t
  = \int_{-\infty}^\infty \theta (x ) \overline{\Omega(x)} d x=: { \la} \theta, \Omega { \ra}
\]
 is correctly defined.

  Note that $h\in  {\cal Z}_{+}'$ if $h\in L^1_{\rm loc}({\Bbb R}_{+})$ and the integral
 \[
\int_{0}^\infty | h(t)| (1+| \ln t |)^{-\kappa } dt< \infty
 \]
 converges  for some $\kappa $. In this case the corresponding function \e{eq:HH3} satisfies the condition
  \[
\int_{-\infty}^\infty |\theta (x)| (1+| x |)^{-\kappa } dx< \infty,
 \]
 and hence $\theta \in{\cal S}'\subset {\cal Z}'$.

    \medskip
 
 {\bf 2.2.}
 Let us now give the definition of the sign-function of a Hankel operator $H$ or of its kernel $h(t)$. Set
  \begin{equation}
 b(\xi) = \frac{1}{ 2\pi }  \frac{ \int_{0}^\infty   h (t) t^{ -i \xi}   dt}{ \int_{0}^\infty   e^{-t} t^{ -i \xi}   dt} .
\label{eq:Bb}\end{equation}
Of course $  b(-\xi)=\ov{ b(\xi)}$ if $  h(t)=\ov{ h(t)}$. 
We call $b(\xi)$ the $b$-function of a Hankel operator $H$ (or of its kernel $h(t)$) and we use the term the {\it  sign-function} for the   Fourier transform $s(x)=\sqrt{2\pi} (\Phi^* b)(x)$ of $b (\xi)$.

Let the function $\theta(\xi)$ be defined by formula  \e{eq:HH3}. If $h \in \cal Z_{+}'$,
then  $\theta \in \cal Z'$ and hence its   Fourier transform 
   \begin{equation}
 a(\xi) =  (\Phi \theta)(\xi)=  (2\pi)^{-1/2}\int_{0}^\infty   h (t) t^{ -i\xi}   dt 
\label{eq:M6}\end{equation}
is      an element  of   $C_{0}^\infty ({\Bbb R})'$.    
Then definition \e{eq:Bb}  can be rewritten  
\begin{equation}
 b(\xi)  =  (2\pi)^{-1/2} a(\xi)  \Gamma (1-i\xi)^{-1}
\label{eq:M9}\end{equation}
where  $\Gamma(\cdot)$ is the gamma function.  
Note that $\Gamma (1-i\xi)^{-1} \neq 0$ for $\xi\in {\Bbb R}$, but according to the Stirling formula it tends exponentially to zero as $|\xi|\to \infty.$  Nevertheless the distribution $b\in C_{0}^\infty ({\Bbb R})'$ and hence $s \in {\cal Z} '$.

For a test function $f \in  {\cal D}$, we set
  \begin{equation}
g(\xi)=\Gamma(1/2+i \xi) (\Phi U f)(\xi)=:( \Xi f)\xi).
\label{eq:MAID1}\end{equation} 
Since $Uf\in{\cal Z}$, the functions $\Phi U f \in C_{0}^\infty ({\Bbb R})$ and hence  $g\in C_{0}^\infty ({\Bbb R})$.
We   note that
  $  (\Phi U f) (\xi)$ is the Mellin transform of $f (t)$.       

The following result was obtained in \cite{Y}.

 \begin{theorem}\label{1}
 Suppose that  $h \in  {\cal Z}_{+}'$. Define the distribution $b\in C_{0}^\infty ({\Bbb R})'$     by formula \e{eq:Bb}, and set $s=\sqrt{2\pi}\Phi^* b\in {\cal Z} '$.
  Let $f_{j}\in {\cal D}$, $j=1,2$,   let the functions      $   g_{j} \in C_{0}^\infty ({\Bbb R})$ be defined by formula \e{eq:MAID1} and $u_{j}=\Phi^* g_{j} =\Phi^* \Xi f_{j} \in {\cal Z} $.  Then      the identity
 \begin{equation}
{ \pmb\la} h ,  \bar{f}_{1}\star f_{2} {\pmb \ra} =   \la s,   \bar{u}_{1} u_2\ra =:  s [u_{1}, u_{2}]   
\label{eq:MAID}\end{equation}
  holds.
 \end{theorem}

       \medskip
 
 {\bf 2.3.}
  For an arbitrary  distribution  $h \in  {\cal Z}_{+}'$, we have constructed in Theorem~\ref{1} its sign-function $s  \in  {\cal Z}'$. 
It turns out that, conversely, the kernel  $h(t)$ can be recovered
from  its sign-function $s(x)$.

    \begin{proposition}\label{round}
  Let    $h \in  {\cal Z}'_{+}$, and let $s \in  {\cal Z}'$ be its sign-function.   Then  
    \begin{equation}
  h (t)  =    \int_{-\infty}^\infty     e^{-t e^{-x}} e^{-x} s (x)  d x.
   \label{eq:con}\end{equation}
  \end{proposition}
        
    As we shall see in the next section, even for kernels \e{eq:FDvm}, the corresponding sign-function $s (x)$ is a highly singular distribution. Nevertheless   the mapping $h(t)\leftrightarrow s (x)$ yields the
        one-to-one correspondence between the classes $ {\cal Z}'_{+}$ and $ {\cal Z}' $.
        We emphasize that formula \e{eq:con} is understood in the sense of distributions.

       \medskip
 
 {\bf 2.4.}
 Suppose  now that $h(t)=\overline{h(t)}$ so that
   the  operator $H$ is self-adjoint.     Then the  identity  \e{eq:MAID}, or equivalently \e{eq:MAIDs}, implies relation \e{eq:MAIDT}.
  To be more precise,  
  we use the following natural definition.  
     Denote by $N_{\pm}(s)$ the maximal dimension of linear sets ${\cal L}_{\pm}\subset {\cal Z}$   such that $\pm s[u,u]  > 0$   for all $u\in {\cal L}_{\pm}$, $u \neq 0$. We apply the same definition to the form $h[f,f]$ considered on the set $  \cal D$ and observe that $N_{\pm}(h)=N_{\pm}(H)$.

    Note that formula \e{eq:MAID1} establishes one-to-one correspondence between the sets $  \cal D$ and $  C_{0}^\infty({\Bbb R})$. Of course  the  Fourier transform establishes one-to-one correspondence between the sets $  C_{0}^\infty({\Bbb R})$ and $  {\cal Z}  $. 
  Therefore the following assertion is a direct consequence of Theorem~\ref{1}.

  \begin{theorem}\label{HBx}
  Let a Hankel operator $H$ be bounded.  Define the distribution $b\in C_{0}^\infty ({\Bbb R})'$     by formula \e{eq:Bb} and set $s=\sqrt{2\pi} \Phi^* b \in  {\cal Z}'$.  Then
     \begin{equation}
  N_{\pm}(H)= N_{\pm}(s).
 \label{eq:Nhs}\end{equation}
   \end{theorem}

In particular, relation \e{eq:Nhs}  means that a Hankel operator $H$ is positive (or negative) if and only if the  function $s(x)$ is positive (or negative). This justifies the term ``sign-function"  for   $s(x)$.

%************************************************************
  %************************************************************
\section{Proofs of Theorems~\ref{FDH1} and \ref{FDHC}}  
 %************************************************************

{\bf 3.1.}
Let us first calculate the $b$- and $s$-functions of the kernel 
       \begin{equation}
h(t)= t^k e^{-\alpha t} \q {\rm where}\q     k =0,1, \ldots, \q \Re\alpha>0,
\label{eq:E1}\end{equation}
but we do not assume that $\Im\alpha=0$.
Calculating integral \e{eq:M6}  we see that
 \[
a(\xi) = (2\pi)^{-1/2} \int_{0}^\infty t^k e^{-\alpha t} t^{-i\xi} dt =
(2\pi)^{-1/2} \alpha^{-1-k + i \xi} \Gamma (1+k-i \xi),
\]
where $\arg \alpha\in (-\pi/2, \pi/2)$,
 and hence  function \e{eq:M9}  equals
 \[
b(\xi) =  \alpha^{-1-k + i \xi} \frac{\Gamma (1+k- i \xi)} { 2\pi  \Gamma (1 -i \xi)}.
\]
Since $k$ is integer,  this yields the following result.

     \begin{lemma}\label{Fr}
     Let $h(t) $  be given by formula \e{eq:E1}. 
      If $k=0$, then $b(\xi) =   (2\pi)^{-1} \alpha^{-1+ i \xi}$ and
        \begin{equation}
s (x) =  \alpha^{-1}\delta (x -\beta), \q \beta=-\ln\alpha.
\label{eq:E6}\end{equation}
 If $k=1,2,\ldots$, then
    \[
b(\xi) =  (2\pi)^{-1 } \alpha^{-1-k+ i \xi} (1-i \xi)\cdots (k-i \xi) 
\]
and
  \begin{equation}
s (x) =  \alpha^{-1-k}(1 -\partial)\cdots (k -\partial)\delta (x -\beta).
\label{eq:E7}\end{equation}
 \end{lemma}

 Let us use the notation $\nu_{\ell,k}$ for the coefficients of the expansion
\[
(1-z)\cdots (k-z)=\sum_{\ell =0}^k \nu_{\ell,k} z^\ell
\]
for $k\geq 1$, $\ell\leq k$,  and set $\nu_{0,0}=1$. Then formulas \e{eq:E6} and \e{eq:E7}  can be rewritten as
\[
s (x) =  \alpha^{-1-k} \sum_{\ell =0}^k \nu_{\ell,k} \delta^{(\ell)} (x -\beta).
\]
Therefore  Lemma~\ref{Fr} implies the following more general result.

 \begin{lemma}\label{FR}
 Let
\begin{equation}
h(t) =  P (t)e^{- \alpha  t}
\label{eq:FDv}\end{equation}
where $\Re \alpha  >0$ and 
  \begin{equation}
   P (t)=\sum_{k=0}^{K }p_{k} t^k  
  \label{eq:Pm}\end{equation}
   is a polynomial. Set
    \begin{equation}
q_{ k} =\sum_{\ell=k }^{K } \nu_{k,\ell} \alpha ^{-1-\ell} p_{\ell} ;
\label{eq:E81q}\end{equation}
 in particular, $q_{K}=(-1)^K \alpha^{-1-K} p_{K}$.
   Then the   $b$- and $s$-functions of   kernel \e{eq:FDv} equal
      \begin{equation}
b(\xi) = (2\pi)^{-1} e^{-ix \beta} Q(x) \q {\rm where} \q Q(\xi)= \sum_{k=0}^{K }q_{ k}   (i\xi)^k,
\label{eq:E81b}\end{equation}
  $ \beta= -\ln \alpha$,  and
   \begin{equation}
s (x) = \sum_{k=0}^{K }q_{ k}   \delta^{(k)} (x -\beta ).
\label{eq:E81v}\end{equation}
  \end{lemma}

 Observe that distribution \e{eq:E81v} is positive
  if and only if $\Im \beta=0$, $q_{ k} =0$ for all $k\geq 1  $ and $q_{ 0}  > 0$. 
 Therefore    the Hankel operator with kernel \e{eq:FDv}, \e{eq:Pm}   cannot be expected to be sign-definite unless $\alpha$  is real and $K =0$.  Theorem~\ref{FDH1} provides essentially more advanced results in this direction. 
 
 Note that for $u\in {\cal Z}$
  \[
    \int_{-\infty}^\infty  \delta^{(k)} (x -\beta ) |u(x)|^2 dx= (-1)^k \sum_{\ell=0}^k C_{k}^\ell     u^{(\ell)}(\beta  )\overline{u^{(k-\ell)}(\bar\beta )}
\]
 where $C_{k}^\ell $ are the binomial coefficients.   This leads to the following result.
   
    \begin{lemma}\label{FR1}
For the distribution given by formula \e{eq:E81v},   we have
  \[
\la s , |u|^2 \ra  
  =  \sum_{j, \ell =0}^{K } s_{j, \ell}  
   u^{(\ell)}(\beta  )\overline{u^{(j)}(\bar\beta )}    ,\q u\in {\cal Z},
\]
  where $s_{j, \ell}    =0$    for $j+\ell   > K $ and
    \begin{equation}
   s_{j, \ell}  = (-1)^{j+\ell } C_{j+\ell }^j q_{j+\ell } 
    \label{eq:rj}\end{equation}
   for $j+\ell  \leq K $; 
  in particular,
   \[
   s_{j, \ell}    = (-1)^{K  }C_{K }^j  q_{K } \q {\rm for}\q j+\ell   = K .
  \]
\end{lemma}

It is now convenient to introduce
   
  \begin{definition}\label{SS}
    Let a kernel $h(t)$ be given by formulas \e{eq:FDv}, \e{eq:Pm}, and let $q_{k}$ be coefficients \e{eq:E81q}. Denote by ${\bf S}  $   the matrix   of order $K+1$  with the elements $ s_{j, \ell}   $ defined in Lemma~\ref{FR1}.  We call ${\bf S} ={\bf S} (P,\alpha)$  the {\it sign-matrix}  of the kernel  $h(t)$. 
\end{definition} 

It is only essential for our proof of Theorem~\ref{FDH1} that the sign-matrix ${\bf S}$ is skew triangular, that is, $s_{j, \ell}    =0$    for $j+\ell   > K $, and that its elements $s_{j, \ell}    = C_{K}^j \alpha^{-1-K} p_{K}$   on the skew-diagonal $j+\ell   = K $ are not zeros if $p_{K}\neq 0$. In this case $\det{\bf S}\neq 0$.   Note  also that 
       \begin{equation}
{\bf S} (\bar{P},\bar\alpha)= {\bf S} (P,\alpha)^*; 
 \label{eq:ZS}\end{equation}
      in particular, $ {\bf S} (P,\alpha)$ is symmetric if $\alpha=\bar{\alpha}$ and $P(t)=\ov{P(t)}$.   
              
          Let us define the mapping $  J_{K} (\beta) : {\cal Z} \to {\Bbb C}^{K +1}$ by the relation\footnote{The upper index ``$\top$" means that a vector is regarded as a column.}
      \begin{equation}
J_{K} (\beta) u= (u (\beta  ),u' (\beta  ), \ldots, u^{(K )}(\beta  ))^\top .
 \label{eq:Z}\end{equation}
 Then   Lemma~\ref{FR1} yields the following assertion.

    \begin{proposition}\label{FR1q}
   For a kernel $h(t)$ defined by \e{eq:FDv}, \e{eq:Pm},   the sign-function is given by the formula
   \begin{equation}
\la s , |u|^2\ra  =  ({\bf S} (P,\alpha)  J_{K} (\beta) u ,   J_{K} (\bar\beta) u )_{K+1} ,\q \beta=-\ln \alpha,  
 \label{eq:Z1}\end{equation}
  where   $( \cdot,  \cdot )_{K+1} $ is   the scalar product  in   ${\Bbb C}^{K +1}$. 
     \end{proposition} 
     
     Formula \e{eq:Z1} is convenient for real $\alpha$ and $P(t)$. In the complex case, we consider the real kernel
 \begin{equation}
h(t) =  P (t)e^{- \alpha  t} +  \ov{P (t)}e^{- \bar\alpha  t}, \q \Re\alpha>0, \q \Im\alpha>0,
\label{eq:FDvc}\end{equation}
corresponding to two complex conjugate points $\alpha$ and $\bar{\alpha}$.  
It follows from Proposition~\ref{FR1q} that the corresponding sign-function equals
\[
\la s , |u|^2\ra  =  ({\bf S} (P,\alpha)  J_{K} (\beta) u ,   J_{K} (\bar\beta) u )_{K+1} 
+ ({\bf S} (\bar{P},\bar\alpha)  J_{K} (\bar\beta ) u ,   J_{K} (\beta) u )_{K+1}.
 \]
 Let us rewrite this equality in the ``matrix" form taking into account relation \e{eq:ZS}.

  \begin{proposition}\label{FR1qc}
   For a kernel $h(t)$ defined by   \e{eq:Pm}, \e{eq:FDvc},   the sign-function is given by the formula 
   \[
\la s , |u|^2\ra   = (  {\bf \wt{S}} (P,\alpha) 
 (  J_{K} (\beta) u , J_{K} (\bar\beta) u )^\top ,
  (  J_{K} (\beta) u , J_{K} (\bar\beta) u )^\top )_{2K+2}   ,    \q \beta=-\ln \alpha, 
\]
       where 
          \begin{equation}
   {\bf \wt{S}} (P,\alpha) = \begin{pmatrix}
0 & {\bf S} (P,\alpha)^*
\\
{\bf S} (P,\alpha)    & 0
\end{pmatrix}.    
 \label{eq:sicc1}\end{equation}
        \end{proposition} 
        
         Let us now consider kernel  \e{eq:FDvm}. We can apply Proposition~\ref{FR1q} to all    real terms  corresponding to $m=1,\ldots, M_{0}$ and Proposition~\ref{FR1qc} to all  complex conjugate terms corresponding to pairs $m$, $  M_1 +m$ where  $m=M_{0}+1,\ldots, M_{1}$. Various objects will be endowed with the index $m=1,\ldots, M_{0}+M_{1}$. Thus we set ${\bf S}_{m}={\bf S} (P_{m}, \alpha_{m})$ for $m=1,\ldots, M_{0}$ and 
         ${\bf S}_{m}={\bf \wt{S}} (P_{m}, \alpha_{m})$ for $m=M_{0}+1,\ldots, M_{0}+M_{1}$.   The mappings ${\bf J}_{m}= J_{K_{m}}(\beta_{m}): {\cal Z}\to {\Bbb C}^{r_{m}}$  are  defined for   $m=1,\dots, M_{0}$ by formula \e{eq:Z} where $ \beta_{m}=-\ln\alpha_{m}$ and $r_{m}=K_{m}+1$.  If  $m=M_{0}+1,\ldots, M_{1}$,  we set
  ${\bf J}_{m}   u= (J_{K_{m}} (\beta_{m})   u, J_{K_{m}} (\bar{\beta}_{m})   u)^\top $; then   ${\bf J}_{m}: {\cal Z}\to {\Bbb C}^{r_{m}}$ where  $r_{m}=2 K_{m}+2$.
 
It is convenient to rewrite the above results in the vectorial notation.
We set
  \begin{equation}
 {\Bbb C}^r=\bigoplus _{m=1}^{M_{0}+   M_1}  {\Bbb C}^{r_m} 
 \label{eq:ccx}\end{equation}
 and introduce the mapping ${\bf J} : {\cal Z}\to  {\Bbb C}^r$ by the formula
 \begin{equation}
 {\bf J} u= ({\bf J}_{1}u,\ldots, {\bf J}_{M_{0}+M_{1}}u)^\top.
 \label{eq:E13bj}\end{equation}
 The sign-matrix of kernel \e{eq:FDvm} is defined as   the block-diagonal matrix
  \begin{equation}
{\bf S}= \diag\{{\bf S}_{1},\ldots, {\bf S}_{M_{0}+M_{1}}\} .
 \label{eq:E13b}\end{equation}
It follows from  Propositions~\ref{FR1q}  and \ref{FR1qc} that the sign-function of kernel  \e{eq:FDvm} is given by the formula
   \begin{equation}
\la s, |u|^2 \ra = ({\bf S}  {\bf J}   u ,  {\bf J}  u )_{r}= \sum_{m=1}^{M_{0}+M_{1}}   ( {\bf S}_{m} {\bf J}_{m}  u ,  {\bf J}_{m} u )_{r_{m}}  .
 \label{eq:si2}\end{equation}

       \medskip         
 
{\bf 3.2.} 
  Below we need the following elementary assertion. We give its proof because  similar arguments will be used  in subs.~3.4 under less trivial circumstances.
  
    \begin{lemma}\label{zz}
    Let $\beta_{1}, \ldots, \beta_{M}\in {\Bbb C}$ and $K_{1},\ldots, K_{M}\in {\Bbb Z}_{+}$. Then there exist functions $\psi_{k, m}\in{\cal Z}$ where $m=1,\ldots, M$ and $k=0,\ldots, K_m$  such that $\psi_{k, m}^{(l)}(\beta_{n}) =\d_{m, n} \d_{k,l}$ for all $n=1,\ldots, M$ and $l=0,\ldots, K_{m}$.
 \end{lemma}
 
   \begin{pf}
   Choose some $m=1, \ldots, M$ and $K \in {\Bbb Z}_{+}$. Let $a_{0},a_{1},\ldots, a_{K}$ be any given numbers.
      It suffices to construct a function $\psi \in {\cal Z}$ such that $\psi ^{(l)}(\beta_{n}) =0$ for all $n\neq m$ and  $\psi ^{(l)}(\beta_m) =a_{l}$ where $l=0,\ldots, K $.

   Let $\varphi_{0} \in {\cal Z}$ be an arbitrary function such that $\varphi_{0}  (0)\neq 0$. Set
   $\omega (z)=1$ if $M=1$, 
   \begin{equation}
\omega (z)=   \prod_{n=1; n\neq m}^M (z-\beta_{n})^{K +1}\q {\rm if }\q M\geq 2,
\label{eq:AF1}\end{equation}
and 
  \begin{equation}
  \varphi (z)=\omega (z) \varphi_{0} (z -\beta_{m}).
\label{eq:AF2}\end{equation}
Of course $\varphi (\beta_{m})\neq 0$. 
  Let us seek the function $\psi $ in the form 
    \begin{equation}
 \psi (z)= Q (z -\beta_{m})\varphi (z)  
\label{eq:AF3}\end{equation}
where   
 \begin{equation}
 Q  (z  )= \sum_{j=0}^{K } q_j z^j
\label{eq:AF4}\end{equation}
 is a polynomial.
Clearly, $\psi \in {\cal Z}$ and $\psi  $ has zeros of   order $K+1$ at all points $\beta_{n}$, $n\neq m$. 

It remains to satisfy the conditions  $\psi ^{(l)}(\beta_m) =a_{l}$. In view of \e{eq:AF3},  \e{eq:AF4} they yield the equations 
\begin{equation}
\sum_{j=0}^l C^j_{l} j! q_j \varphi ^{(l-j)}(\beta_m) =a_{l}, \q l=0, 1,\ldots, K , 
\label{eq:AF5}\end{equation}
for the coefficients $q_j$.
For $l =0$,  we find that 
\begin{equation}
q_{0 } = \varphi (\beta_m)^{-1} a_{0}  .
\label{eq:AF5w}\end{equation}
Then equation \e{eq:AF5} determines $q_{l }$ if $q_{0 },\ldots, q_{l-1 }$ are already found. The corresponding function \e{eq:AF3} satisfies all necessary conditions. 
   \end{pf}
   
    Set    ${\sf u}_{m}= (u_{0,m}, u_{1,m}, \ldots, u_{K_{m},m})^\top \in {\Bbb C}^{K_{m}+1}$ for $m= 1,\ldots, M $,
   ${\bf u}_{m}={\sf u}_{m}$ for $m= 1,\ldots, M_{0} $ and  ${\bf u}_{m}=({\sf u}_{m}, {\sf u}_{m+M_{1}})^\top$ for $m= M_{0}+ 1,\ldots, M_{0} +M_{1}$. 
    Then  ${\bf u}_{m} \in {\Bbb C}^{r_m}$   and ${\bf u}=({\bf u}_{1}, \ldots, {\bf u}_{M_{0} +M_{1}})^\top$ is an element of the direct sum  \e{eq:ccx}. Let us define the mapping ${\bf Y}: {\Bbb C}^r\to {\cal Z}$ by the formula
 \begin{equation}
 ({\bf Y} {\bf u})(z)= \sum_{m=1}^M \sum_{k=0}^{K_m} u_{k,m}\psi_{k, m} (z)   
\label{eq:YY}\end{equation}
where $\psi_{k, m}$ are the functions constructed in Lemma~\ref{zz}.
We apply this definition in the case where $\beta_m=\bar{\beta}_m$ for $m=1,\ldots, M_0$ and
$\beta_m=\bar{\beta}_{m+M_1}$, $K_m=K_{m+M_1}$ for $m=M_0+1,\ldots, M_0+M_1$. 
By the definition of the functions $\psi_{k, m} $, for mapping \e{eq:E13bj} we have the identity
 \begin{equation}
 {\bf J} {\bf Y} =I.   
\label{eq:YY1}\end{equation}

In view of Theorem~\ref{HBx}, for the proof of Theorem~\ref{FDH1} we   only have to calculate the numbers $ N_{\pm}(s)$. This can be reduced to  a problem of the linear algebra.

 \begin{lemma}\label{yy}
  Let $s$ be the sign-function of    kernel  \e{eq:FDvm}, and let $\bf S$ be the corresponding sign-matrix defined by formula \e{eq:E13b}. Then
    \begin{equation}
 N_{\pm}(s)=  N_{\pm}({\bf S}).  
\label{eq:YY2}\end{equation}
 \end{lemma}
 
  \begin{pf}
We proceed from identity  \e{eq:si2}.  Consider, for example, the sign $``-"$.  
        If  $\la s,  |u|^2\ra<0$,   then    $({\bf S}{\bf u},  {\bf u})_{r}<0$ for ${\bf u}= {\bf J} u$.  This shows that $N_{-} (s)\leq N_{-}({\bf S})$.
 
 Let us prove the opposite inequality.     It follows from the identities \e{eq:si2} and  \e{eq:YY1} that
 \[
 \la s, |{\bf Y} {\bf u}|^2 \ra = ({\bf S}  {\bf u} ,  {\bf u} )_{r} .
 \]
 Thus if  $({\bf S}{\bf u},  {\bf u})_{r}<0$,   then   $\la s,  |u|^2\ra<0$ for
  $u= {\bf Y} {\bf u}$.  
    \end{pf}

  \medskip
   
  {\bf 3.3.}
  It remains to calculate the numbers
   \begin{equation}
 N_{\pm}({\bf S})=\sum_{m=1}^{M_{0}+M_{1}}  N_{\pm}({\bf S}_{m}).  
\label{eq:YY3}\end{equation}
It is quite easy to find $N_{\pm}({\bf S}_{m})$ for $m\geq M_{0} + 1  $.

  \begin{lemma}\label{SThc}
  Under the assumptions of Proposition~\ref{FR1qc} suppose that $p_{K}\neq 0$.
 Then matrix \e{eq:sicc1}   has exactly $K+1$ positive and $K+1$ negative eigenvalues $($they are opposite to each other$)$.  
 \end{lemma}  
 
   \begin{pf} 
   Set  ${\bf S}= {\bf S} (P,\alpha)$ and  recall that $\det {\bf S} \neq 0$. If ${\bf S}^* {\bf S}f = \lambda^2 f $ for some $\lambda>0$, then
\[
\wt{\bf S} \begin{pmatrix}
\lambda f 
\\
\pm {\bf S}    f 
\end{pmatrix}=
 \begin{pmatrix}
0 & {\bf S} ^*
\\
{\bf S}     & 0
\end{pmatrix}  
\begin{pmatrix}
\lambda f 
\\
\pm {\bf S}    f 
\end{pmatrix}=\pm \lambda \begin{pmatrix}
\lambda f 
\\
\pm {\bf S}    f  \end{pmatrix} .
\]
Thus we put into correspondence to every eigenvalue $\lambda^2$ of the matrix  ${\bf S}^* {\bf S}   $ of order $K+1$ the  eigenvalues $\lambda $ and $-\lambda $ of the matrix $\wt{\bf S}$ of order $2K+2$.
\end{pf}

  In the case  $m \leq M_{0}$ we need some information on skew triangular matrices. We consider Hermitian matrices $S$ of order $K+1$ with elements $s_{j,\ell}$, $ j,\ell= 0,\ldots,K$, such that $s_{j,\ell}=\bar{s}_{\ell,j}$.
  We say that a matrix $S$  is skew triangular if $s_{j,\ell}=0$ for $ j+\ell> K$. It is easy to see (reasoning, for example, by induction) that
  \begin{equation}
\det S= (-1)^{K (K+1)/2} s_{0,K} s_{1,K-1}\cdots s_{K,0}.
\label{eq:ST}\end{equation}
In particular, $\det S \neq 0$ if (and only if) all skew diagonal elements are not zeros.

Let us first consider    skew diagonal matrices.

 \begin{lemma}\label{STa}
Let $S_{0}$ be a   Hermitian    matrix of order $K+1$ such that   $s_{j,\ell}=0$ for  $ j +\ell\neq K $.  If $K$ is odd, then $S_{0}$ has the eigenvalues $\pm |s_{j,K-j}|$ where $j=0,\ldots, (K-1)/2$.  If $K$ is even, then $S_{0}$ has 
the eigenvalues $\pm |s_{j,K-j}|$ where $j=0,\ldots, K /2-1$ and the eigenvalue $s_{K/2, K/2}$. 
 \end{lemma}
 
 \begin{pf}
 Let us consider the equation $S_{0}f=\lambda f$ for $f=(f_{0},\ldots, f_K)^\top$. 
Since  $S_{0} f=(s_{0,K} f_K, s_{1,K-1} f_{K-1}, \ldots, s_{ K,0}f_0)^\top$   this equation   is equivalent to the system
  \begin{equation}
s_{j, K-j} f_{K-j}=\lambda f_{j},\q j=0,\ldots, K.
\label{eq:ST1}\end{equation}
If $K$ is odd, then \e{eq:ST1} decouples into  $(K+1)/2$ systems of two equations   for $f_{j}$ and $f_{K-j}$ where $j=0,\ldots, (K-1)/2$. Every such system has two simple eigenvalues $\lambda= \pm  \sqrt{s_{j,K-j} s_{K-j,j}} = \pm |s_{j,K-j}|$.
If $K$ is even, then \e{eq:ST1} decouples into  $K/2$ systems of the same two equations   for $f_{j}$ and $f_{K-j}$ where $j=0,\ldots, K/2-1$ and the single equation $s_{K/2,K/2}  f_{K/2}=\lambda f_{K/2}$. The last equation has of course the eigenvalue $\lambda= s_{K/2,K/2}$. 
   \end{pf}
 
For applications to Hankel operators, we need the following result.

  \begin{lemma}\label{ST}
Let $S$ be a Hermitian     skew triangular matrix of order $K+1$ such that $s_{j, K-j}\neq 0$ for $j=0,\ldots, K$. If $K$ is odd, then $S$ has $(K+1)/2$ positive and $(K+1) /2$ negative eigenvalues. If $K$ is even, then $S$ has $K/2+1$ positive and $K/2$ negative eigenvalues for $s_{K/2,K/2}>0$ and it has $K/2 $ positive and $K/2+1$ negative eigenvalues for $s_{K/2,K/2}<0$.
 \end{lemma}
 
 \begin{pf}
 According to formula \e{eq:ST},  $\det S$ depends only on elements  $s_{j,\ell}$ on the skew diagonal where  $ j +\ell= K$. Let us use that eigenvalues of $S$ depend continuously on its matrix elements  so that they cannot cross the point zero unless one of  skew diagonal elements hits the  zero.
 
  Let us consider the family of matrices $S(\varepsilon)$ where $\varepsilon\in [0,1]$  with  elements $s_{j,\ell}(\varepsilon)=\varepsilon s_{j,\ell}$ for $ j +\ell < K $ and  $s_{j,\ell}(\varepsilon)=  s_{j,\ell}$ for $ j +\ell \geq K $. Since $\det S(\varepsilon)= \det S\neq 0$ for $\varepsilon\in [0,1]$,
  all matrices $S(\varepsilon)$ and, in particular, $S(1)=S$ and $S(0)$, have the same numbers of positive and negative eigenvalues.   So it remains to apply  Lemma~\ref{STa} to the matrix $S(0)$. 
   \end{pf}

    The following result is a particular case of Lemma~\ref{ST}.

  \begin{lemma}\label{STh}
  Let     ${\bf S} ={\bf S} (P,\alpha)$ be the sign-matrix of kernel \e{eq:FDv}, \e{eq:Pm} where $\Im\alpha=0$,  $ P(t) = \ov{P(t)}$  and $p_{K}\neq 0$.
    The total numbers ${\cal N}_{+} =N_{+}({\bf S})$ and ${\cal N}_{-} =N_{-}({\bf S})$ of strictly positive   and negative    eigenvalues of the matrix ${\bf S} $ are given by the equalities
\begin{align*}
{\cal N}_{+}  = {\cal N}_{-} &=(K +1)/2 \q {\rm if} \q K \q {\rm is  } \; {\rm   odd}
 \\
{\cal N}_{+}  -1= {\cal N}_{-}  &= K  /2 \q {\rm if} \q K \q {\rm is  } \; {\rm   even} \; {\rm   and}
 \q p_{K }  > 0
 \\
{\cal N}_{+}= {\cal N}_{-}  -1  &= K  /2  \q {\rm if} \q K \q {\rm is  } \; {\rm   even} \; {\rm   and}
 \q p_{K}  < 0 .
 \end{align*}
 \end{lemma}
  
Combined with equality \e{eq:YY3},
Lemmas~\ref{SThc}   and \ref{STh}  show that     
 \begin{equation}
 N_{\pm} ({\bf S}) =  \sum_{m=1}^{M_{0}}  {\cal N}_\pm^{(m)} +  \sum_{m=M_{0} + 1}^{M_{0}+M_{1}}  K_{m} + M_{1} .
   \label{eq:TNL}\end{equation}
   Putting  this result  together  with relations \e{eq:Nhs}  and \e{eq:YY2}, we
       conclude the proof of Theorem~\ref{FDH1}. 

  \medskip 
   
 {\bf 3.4.}
In this subsection we consider operators $H=H_{0}+V$ where $H_{0}$
 is the Carleman operator (or a more general operator) and   $V$ is a finite rank Hankel operator. Various objects related to the operator $H_{0}$ will be endowed with the index $``0"$,  and objects related to the operator $V$ will be endowed with the index $``v"$. Our goal is 
get an explicit formula for the total number $N_- (H)$ of negative eigenvalues of the operator $H $. 
 
  \begin{theorem}\label{FDH2}
 Suppose that the sign-function $s_{0} (x)$ of a Hankel operator $ H_{0}$ is    bounded and positive.   Let the kernel $v(t)$ of $V$ be given by the formula
  \[
v(t) =\sum_{m=1}^M  P_{m} (t)e^{- \alpha_{m}  t}
\]
where $P_{m} (t)$ is a polynomial of degree $K_{m}$.
 Define the numbers ${\cal N}_- ^{(m)}$   by formula \e{eq:RXm} where ${\sf p}_{m}$ is   coefficient \e{eq:ppq}. 
  Then  the total number $N_- (H)$ of negative eigenvalues of the operator   $H=H_{0}+V$      is given by formula  \e{eq:TN}.
   \end{theorem} 
 
Comparing Theorem~\ref{FDH1}  for the operator $V$ and Theorem~\ref{FDH2}, we 
can state the following result.

   \begin{theorem}\label{FDH}
   Under the assumptions of  Theorem~\ref{FDH2}, we have
   $$N_{-} (H)=N_{-} (V).$$
    In particular, $H \geq 0$ if and only if $V\geq 0$.
 \end{theorem}

 Since for the Carleman operator ${\bf C}$ the sign-function  $s_{0}(x)=1$, Theorem~\ref{FDH} applies to $H_{0}={\bf C}$ and hence implies Theorem~\ref{FDHC}.

  The proof of Theorem~\ref{FDH2} is essentially similar to that of Theorem~\ref{FDH1}.
  Relation \e{eq:Nhs} remains of course true but instead of \e{eq:si2} we now have
     \begin{equation}
\la s, |u|^2 \ra = \int_{-\infty}^\infty s_{0}(x)Ê|u(x)|^2 dx +  ({\bf S}_{v}  {\bf J}   u ,  {\bf J}  u )_{r}  .
 \label{eq:E12z}\end{equation}
  Compared to subs.~3.2,  we additionally  have   to     consider
  the first term in the right-hand side of  \e{eq:E12z}. Instead of Lemma~\ref{zz}, this requires a more special assertion.
  
   \begin{lemma}\label{zzb}
    Let $\beta_{1}, \ldots, \beta_{M}\in {\Bbb C}$,  $K_1,\ldots, K_{M}\in {\Bbb Z}_{+}$ and $\varepsilon> 0$. Then there exist functions $\psi_{k, m}(\varepsilon)\in{\cal Z}$ where $m=1,\ldots, M$ and $k=0,\ldots, K_m$  such that $\psi_{k, m}^{(l)}(\beta_{n}; \varepsilon) =\d_{m, n} \d_{k,l}$ for all $n=1,\ldots, M$ and $l=0,\ldots, K_m$. Moreover, these functions satisfy the condition
         \begin{equation}
 \int_{-\infty}^\infty  |\psi_{k, m}(x; \varepsilon)|^2 dx = O (\varepsilon),\q \varepsilon\to 0 .
 \label{eq:YY6}\end{equation}
 \end{lemma}
  
 \begin{pf}
    Choose some $m=1, \ldots, M$ and $K \in {\Bbb Z}_{+}$. Let $a_{0},a_{1},\ldots, a_{K}$ be any given numbers.
      It suffices to construct a function $\psi (\varepsilon)\in {\cal Z}$ such that $\psi ^{(l)}(\beta_{n}; \varepsilon) =0$ for all $n\neq m$ and  $\psi ^{(l)}(\beta_m; \varepsilon) =a_{l}$ where $l=0,\ldots, K $.  We  also have  to satisfy condition \e{eq:YY6} for the function $\psi (x; \varepsilon)$.
     
      Instead of \e{eq:AF2} we now set
       \begin{equation}
  \varphi (z;\varepsilon)= \omega (z) \varphi_{0} ((z -\beta_{m})/\varepsilon)
\label{eq:AF2b}\end{equation}
where, as in  Lemma~\ref{zz},   $\varphi_{0}  \in {\cal Z}$ is an arbitrary function such that $\varphi_{0} (0)\neq 0$ and $\omega (z)$ is function \e{eq:AF1}.
We again  seek the function $\psi (\varepsilon)$ in the form 
    \begin{equation}
 \psi (z;\varepsilon)= Q (z -\beta_{m};\varepsilon)\varphi (z;\varepsilon)  
\label{eq:AF3b}\end{equation}
where  $Q $ is polynomial \e{eq:AF4} with the coefficients $q_j =q_j (\varepsilon)$   depending on $\varepsilon$.
As before, $\psi (\varepsilon)\in {\cal Z}$ and $\psi (\varepsilon) $ has zeros of  order $K+1$ at all points $z_{n}$ for $n\neq m$ and all $\varepsilon>0$. 

The conditions  $\psi ^{(l)}(\beta_m,\varepsilon) =a_l$ yield again  equations \e{eq:AF5},  but now the coefficients   $ \varphi ^{(l-j)}(\beta_m; \varepsilon)$ depend on $\varepsilon$. Note that $ \varphi  (\beta_m; \varepsilon)= \omega  (\beta_m) \varphi_{0} (0) $ and according to \e{eq:AF2b}
\begin{equation}
|\varphi ^{(k)} (\beta_m; \varepsilon)|Ê\leq C_k\varepsilon^{-k}.
\label{eq:AF5b}\end{equation}
The coefficient $q_{0}$ is again determined by formula \e{eq:AF5w}; it does not depend on $\varepsilon$. Solving equations \e{eq:AF5} successively for $q_{1} ( \varepsilon) ,\ldots, q_{K} ( \varepsilon)$ and using estimates  \e{eq:AF5b} we find that
\begin{equation}
|q_j ( \varepsilon)|Ê\leq C_j \varepsilon^{-j},\q j=0,\ldots, K.
\label{eq:AG5b}\end{equation}
It follows from \e{eq:AF4}, \e{eq:AF2b} and  \e{eq:AF3b}  that, for   $N= (K+1) (M-1)$, 
\begin{multline*}
  \int_{-\infty}^\infty  |\psi (x; \varepsilon)|^2 dx
  \\
   \leq C \sum_{j=0}^{K}|q_j (\varepsilon)|^2 \int_{-\infty}^\infty (x-\beta_{m})^{2j} (1+ |x- \beta_{m}|^{2N}) |\varphi_{0} ((x-\beta_{m})/\varepsilon)|^2
 dx .
  \end{multline*}
  The integrals in the right-hand side are bounded by $C \varepsilon^{2 j+1}$.
 In view of \e{eq:AG5b} it follows that this expression is $O(\varepsilon)$ as $\varepsilon\to 0$.
   \end{pf}

    Let us return to  Theorem~\ref{FDH2}. Recall that the sign-function of the operator $H$ is given by equality  \e{eq:E12z}.
        According to Theorem~\ref{HBx} and formula \e{eq:TNL} for $N_{-}({\bf S}_{v})$,  we only have to check that $N_{-} (s) = N_{-}({\bf S}_{v})$. 
         Since $s_{0}(x)\geq 0$, we have the estimate
          \[
\la s, |u|^2 \ra \geq  ({\bf S}_{v}  {\bf J}   u ,  {\bf J}  u )_{r}  
 \]
  which directly implies (cf. the proof of Lemma~\ref{yy})    that $N_{-} (s)\leq N_{-}({\bf S}_{v})$.
         
          It remains to check that     $N_{-} (s)\geq N_{-}({\bf S}_{v})$. Let $\bf L$ be the subspace of ${\Bbb C}^r$ spanned by the eigenvectors of ${\bf S}_{v}$  corresponding to its negative eigenvalues. Then $\dim{\bf L}=N_{-}({\bf S}_{v})$ and there exists $\lambda_{0}>0$ such that
    \begin{equation}
     ({\bf S}_{v}{\bf u},{\bf u})_{r}\leq -\lambda_{0} \| {\bf u} \|_{r}, \q \forall \bf u\in \bf L.
   \label{eq:BG5b}\end{equation}
          We again define the function $u (\varepsilon)= {\bf Y}  (\varepsilon){\bf u}$ by formula \e{eq:YY}  where $\psi_{k,m}(z,\varepsilon)$ are the functions constructed in Lemma~\ref{zzb} for sufficiently small $\varepsilon $. Similarly to  \e{eq:YY1}, we have ${\bf J}  {\bf Y}  (\varepsilon)=I$.          Since $s_0\in L^\infty ({\Bbb R})$, it follows from equality  \e{eq:E12z} and estimates  \e{eq:YY6},  \e{eq:BG5b} that
          \[
      \la s,  |u (\varepsilon) |^2\ra\leq -(\lambda_{0} - C\varepsilon ) \| {\bf u} \|_{r}, \q \forall \bf u\in \bf L.
        \]   
     Choosing $C\varepsilon<\lambda_{0}$, we see that $N_{-} (s)\geq \dim{\bf L}$.  This concludes the proof of Theorem~\ref{FDH2}.
 
  \medskip  

 {\bf 3.5.}
 It follows from Lemma~\ref{FR}  that  the $b$- and $s$-functions of  kernel  \e{eq:FDvm} are   the sums (over $m$) of terms \e{eq:E81b} and \e{eq:E81v}, respectively. The coefficients $q_{k,m}$ of the corresponding polynomials $Q_{m}(\xi)$ are constructed  by formula \e{eq:E81q} in terms of the coefficients $p_{k,m}$ of the   polynomials $P_{m}(t)$.

 It turns out that formulas \e{eq:E81b} or  \e{eq:E81v}   for the $b$- or  $s$-functions characterize finite rank Hankel operators. Moreover, the coefficients of the polynomials $P_{m}(t)$ are
 determined by   the coefficients of the polynomials $Q_{m}(\xi)$. This follows from the assertion below.
 
   \begin{lemma}\label{STrb}
   If a function $b(\xi)$ is defined by formula \e{eq:E81b}, then  there exists the unique polynomial $P (t)$ of degree $K$    such that $b(\xi)$   is the $b$-function of the kernel $h (t)=P  (t)e^{-\alpha t}$ where $\alpha=e^{-\beta}$.
 \end{lemma}
 
  \begin{pf}
  Let us  solve equations \e{eq:E81q} for the coefficients $p_{k} $ where $k=K,K-1,\ldots, 0$. Recall  that $\nu_{k,k}=(-1)^k$. Therefore   according to  equation \e{eq:E81q} where $k=K $ we have  
  \[
  p_{K } =(-1)^{K }  \alpha^{1+K } q_{K } 
  \]
and, more generally, 
     \[
p_{ k} = (-1)^K \alpha^{1+K} q_{k} -  (-1)^K  \sum_{\ell=k+1 }^{K } \nu_{k,\ell} \alpha ^{k-\ell} p_{\ell}.
\]
Thus we can   successively find all coefficients $p_{K}, p_{K-1}, \ldots, p_{0}$. 
    \end{pf}
    
    Since the functions $b(\xi)$ and $s(x)$ are obtained from each other by the Fourier transform, Lemma~\ref{STrb} can be equivalently reformulated in terms of the sign-functions.   Of course the reconstructions of $h(t)$ by formula \e{eq:con} and by the method of Lemma~\ref{STrb} are consistent with each other.

    Finally, we   state an equivalent assertion in terms of the sign-matrices $\bf S$ (see Definition~\ref{SS}). We recall that the sign-matrix $\bf S=\bf S (P,\alpha)$ of the kernel $h(t)$ defined by \e{eq:FDv}, \e{eq:Pm}  is skew triangular and  its matrix elements $s_{j, \ell}    = C_{K}^j \alpha^{-1-K} p_{K}$   if $j+\ell   = K $. As usual, we suppose that $p_{K}\neq 0$. 
  Moreover, the matrix ${\bf S} $ possesses an additional property: the numbers 
  \begin{equation}
   \ell ! j! s_{\ell, j} = :\rho_{\ell + j},\q  \ell, j=0,1,\ldots, K,
\label{eq:rr}\end{equation}
  depend on the sum $\ell+ j$ only. We also note that the matrix $\bf S$ is symmetric if $ \alpha =\bar{\alpha}$ and $P(t)=\ov{P(t)}  $.

  The following assertion shows that there is one-to-one correspondence between Hankel kernels $h (t)=P  (t)e^{-\alpha t}$ and such matrices.
  
   \begin{lemma}\label{STr}
 Let   elements $s_{\ell, j}  $ of  a  skew triangular matrix $S $ of order $K+1$ satisfy condition \e{eq:rr}. Then, for every $\alpha$ with $\Re\alpha>0$,  there exists the unique polynomial $P (t)$ of degree $K$    such that $S= {\bf S}(P,\alpha) $ is the sign-matrix of the kernel $h (t)=P  (t)e^{-\alpha t}$.
 \end{lemma}
 
  \begin{pf}
  Comparing relations \e{eq:rj} and \e{eq:rr}, we see that   $q_{k} = (-1)^k k!^{-1}\rho_{k} $. Thus it remains  to use Lemma~\ref{STrb}. 
    \end{pf}

%%%%%%%%%%%%%%%%%%%%%%%%%%%%%%%%%%%%%
 \section{Various  representations of Hankel operators}
 %%%%%%%%%%%%%%%%%%%%%%%%%%%%%%%%%%%%

 Hankel operators can be realized in various spaces. We distinguish four representations: in the spaces ${\Bbb H}_{+}^2 ({\Bbb T})$,    
${\Bbb H}_{+}^2 ({\Bbb R})$,  $l^2  ({\Bbb Z}_{+})$  and $L^2 ({\Bbb R}_{+})$. The last one was already used above.
For the precise definitions of the Hardy classes ${\Bbb H}_{+}^2 ({\Bbb T})$ and    
${\Bbb H}_{+}^2 ({\Bbb R})$, see, e.g., the book \cite{Hof}.
In this text we describe   bounded  Hankel operators  in terms of their quadratic forms. Our presentation seems to be somewhat different from those in the books \cite{Pe,Po}.

  \medskip

 {\bf 4.1.}
 Let us start with the representation of Hankel operators in the  Hardy space 
${\Bbb H}^2_{+} ({\Bbb T})\subset L^2 ({\Bbb T})$ 
of functions analytic in the unit disc. An operator $\bf G$ in the space ${\Bbb H}^2_{+} ({\Bbb T})$ is called Hankel if its quadratic form admits the representation
\begin{equation}
 ( {\bf G } u,u)= \int_{\Bbb T} \omega (\mu) u(\bar{\mu}) \ov{u(\mu)}    d{\bf m}(\mu),\q d{\bf m}(\mu)=(2\pi i \mu)^{-1}d\mu, \q \forall u\in {\Bbb H}^2_{+} ({\Bbb T}) ,
\label{eq:QF}\end{equation}
where $\omega\in L^\infty ({\Bbb T})$.
Note that $d{\bf m}(\mu)$ is 
the Lebesgue measure on ${\Bbb T}$ normalized so that ${\bf m} ({\Bbb T})=1$. 
The operator $ {\bf G } $ is determined by the  function $\omega(\mu)$, 
that is, $ {\bf G } = {\bf G }  (\omega)$. 
The function $\omega(\mu)$ is known as the symbol of the Hankel operator
 $  {\bf G}  (\omega)$. 
 Of course the symbol is not unique because $  {\bf G }  (\omega_{1})
 =  {\bf G }  (\omega_2)$ if (and only if) $\omega_{1}-\omega_{2}\in {\Bbb H}^\infty_{-} ({\Bbb T})$ (the space of  analytic functions outside of the unit disc  bounded and decaying at infinity).
 
  Hankel operators in the  Hardy space 
${\Bbb H}^2_{+} ({\Bbb R})\subset L^2 ({\Bbb R})$ 
of functions analytic in the upper half-plane are defined quite similarly. An operator $\bf H$ in the space ${\Bbb H}^2_{+} ({\Bbb R})$ is called Hankel if its quadratic form admits the representation
\begin{equation}
 ( {\bf H } w,w)= \int_{\Bbb R} \varphi (\lambda) w(-\lambda) \ov{w(\lambda)}    d\lambda, \q \forall w\in {\Bbb H}^2_{+} ({\Bbb R}) ,
\label{eq:QFR}\end{equation}
where $\varphi\in L^\infty ({\Bbb R})$.  
The operator $ {\bf H} $ is determined by the  function $\varphi(\lambda)$, 
that is, $ {\bf H } = {\bf H}  (\varphi)$. 
The function $\varphi(\lambda)$ is known as the symbol of the Hankel operator
 $  {\bf H }  (\varphi)$. 
 Of course the symbol is not unique because $  {\bf H }  (\varphi_{1})
 =  {\bf H }  (\varphi_2)$ if (and only if) $\varphi_{1}-\varphi_{2}\in {\Bbb H}^\infty_{-} ({\Bbb R})$ (the space of bounded analytic functions in the lower half-plane).
 
 In the   space 
$l^2 ({\Bbb Z}_{+})$ of sequences $\xi=(\xi_{0},\xi_{1},\ldots)$,  a   Hankel operator $G$ is defined via its quadratic form
\begin{equation}
 ( G \xi,\xi)= \sum_{n,m=0}^\infty \varkappa_{n+m}\xi_{m}\bar{\xi}_{n}.
\label{eq:QFK}\end{equation}
It is first considered on vectors $\xi$ with only a finite number of non-zero components, and it is supposed that 
\begin{equation}
\big| \sum_{n,m=0}^\infty \varkappa_{n+m}\xi_{m}\bar{\xi}_{n} \big| \leq C \| \xi \|^2.
\label{eq:QFK1}\end{equation}
  Then there exists a bounded operator $G$ such that relation \e{eq:QFK}  holds. We note that condition \e{eq:QFK1} directly implies that $\varkappa=(\varkappa_{0}, \varkappa_{1},\ldots)\in l^2 ({\Bbb Z}_{+})$. Indeed, passing from the quadratic form for $\xi$ to the sesquilinear form for $\xi,\eta$ and choosing $\eta=(1, 0,0, \ldots)$, we see that
  \[
  \big| \sum_{n=0}^\infty \varkappa_{n }\xi_{m}  \big| \leq C \| \xi \|,
  \]
whence  $\varkappa \in l^2 ({\Bbb Z}_{+})$.

Finally  in the   space 
$L^2 ({\Bbb R}_{+})$,    a   Hankel operator $H$ is defined via its quadratic form
\begin{equation}
 ( H f, f)=  {\pmb\la} h, \bar{f}\star f {\pmb\ra}
\label{eq:QFH}\end{equation}
where $f\in C_{0}^\infty ({\Bbb R}_{+})$, $ \bar{f}\star f $ is the Laplace convolution \e{eq:HH1} (note that  $ \bar{f}\star f \in C_{0}^\infty ({\Bbb R}_{+})$) and the distribution $h\in C_{0}^\infty ({\Bbb R}_{+})'$. If
\begin{equation}
|   {\pmb\la} h, \bar{f}\star f {\pmb\ra}|\leq C \| f \|^2,
\label{eq:QFH1}\end{equation}
then there exists a bounded operator $H$ such that relation \e{eq:QFH}  holds.  

It is easy to see that Hankel operators $\bf G$, $\bf H$, $G$ and $H$ are self-adjoint  if  
          \[
    \omega(\bar{\mu} )=\overline{\omega(\mu)}, \q \varphi (-\lambda)=\overline{\varphi(\lambda)}, \q  \varkappa_{n}=\bar{\varkappa}_{n}   \q{\rm and }\q  h(t ) = \overline{h (t )},
    \]
 respectively.

 \medskip

 {\bf 4.2.}
 Let us establish one-to-one correspondences between the  representations of Hankel operators in the spaces ${\Bbb H}_{+}^2 ({\Bbb T})$,    
${\Bbb H}_{+}^2 ({\Bbb R})$,  $l^2  ({\Bbb Z}_{+})$  and $L^2 ({\Bbb R}_{+})$. Let us introduce the notation ${\Bbb A}({\cal H})$ for the linear set of all bounded Hankel operators acting in one of these four Hilbert spaces $\cal H$. It is easy to see that  $H^*\in{\Bbb A}({\cal H})$ together with $H$ and that ${\Bbb A}({\cal H})$ is a closed set in the weak operator topology.
 
 Recall that the function
\begin{equation}
\zeta=\frac{z -i}{ z +i}
\label{eq:chv}\end{equation}
determines a conformal mapping $z\mapsto\zeta$ of the upper half-plane onto the unit disc. The unitary operator  ${\cal U}: {\Bbb H}_{+}^2({\Bbb T})\to  {\Bbb H}_{+}^2 ({\Bbb R} )$ corresponding to this mapping is defined by the equality
\begin{equation}
({\cal U} u)(\lambda)= \pi^{-1/2} (\lambda+i)^{-1} u \bigl(\tfrac{\lambda-i}{\lambda+i}\bigr).
\label{eq:pwtq}\end{equation}
Making the change of variables  \e{eq:chv} in \e{eq:QF}, we see that 
 \[
{\bf G}\in {\Bbb A}( {\Bbb H}_{+}^2 ({\Bbb T})) \Longleftrightarrow {\bf H} ={\cal U}{\bf G}{\cal U}^*\in {\Bbb A}( {\Bbb H}_{+}^2 ({\Bbb R}))
\]
  if   the symbols $ \omega$ of $\bf G$ and $\varphi$ of ${\bf H}$  are linked by the formula
\begin{equation}
\varphi(\lambda)=-\frac{\lambda-i}{\lambda +i}\: \omega\big(\tfrac{\lambda -i}{\lambda+i}\big).
 \label{eq:pw}\end{equation}

The unitary mapping ${\cal F}: {\Bbb H}_{+}^2 ({\Bbb T})\to l^2({\Bbb Z}_{+})$ corresponds 
to expanding  a function in the Fourier series:
\[
\xi_{n}=({\cal F} u)_{n} =  \int_{\Bbb T} u(\mu)    \mu^{-n } d {\bf m}(\mu).
\]
Conversely, for a sequence $\xi= \{\xi_{n}\}$, we have
\[
u(\mu) = ({\cal F}^* \xi) (\mu) = \sum_{n=0}^\infty \xi_{n}\mu^n.
\]
Substituting this expansion into \e{eq:QF}, we see that
\begin{equation}
 ( {\bf G } u,u)=  \sum_{n,m=0}^\infty \varkappa_{n+m}\xi_{m}\bar{\xi}_{n}
 \label{eq:CF1p}\end{equation}
  where 
  \begin{equation}
\varkappa_{n}=({\cal F} \omega)_{n}, \q n\in {\Bbb Z}_{+}.
\label{eq:pww}\end{equation}
If $\omega\in L^\infty ({\Bbb T})$, then expression  \e{eq:CF1p} satisfies estimate \e{eq:QFK1},  and hence the operator $G$ defined by relation \e{eq:QFK} is bounded.   Thus the inclusion $ {\bf G} \in {\Bbb A}( {\Bbb H}_{+}^2 ({\Bbb T}))$ implies that
 $G= {\cal F}{\bf G}{\cal F}^* \in {\Bbb A}( l^2 ({\Bbb Z}_{+}))$. 
 
  Conversely, according to the   Nehari theorem (see the original paper \cite{Nehari}, or the books 
   \cite{Pe}, Chapter~1, \S 1 or \cite{Po}, Chapter~1, \S 2) under assumption \e{eq:QFK1} there exists a function $\omega\in L^\infty ({\Bbb T})$ such that equalities \e{eq:pww} hold. Hence the operator $\bf G$ defined  by relation  \e{eq:CF1p} satisfies also \e{eq:QF}.   Thus for every $ G  \in {\Bbb A}( l^2 ({\Bbb Z}_{+}))$ the operator  ${\bf G} = {\cal F}^*  G{\cal F} \in  {\Bbb A}( {\Bbb H}_{+}^2 ({\Bbb T}))$. 
   
   To show that the inclusions $ {\bf H} \in {\Bbb A}( {\Bbb H}_{+}^2 ({\Bbb R}))$ and 
    $H= \Phi {\bf H} \Phi^* \in {\Bbb A}( L^2 ({\Bbb R}_{+}))$ are equivalent, we need the continuous version of the   Nehari theorem. We give only its brief proof referring to subs.~3.2 of \cite{Y} for details.
    
      \begin{proposition}\label{NehC}
      Let $h\in C_{0}^\infty ({\Bbb R}_{+})'$. Then estimate  \e{eq:QFH1} is satisfied if and only if there exists a function $\varphi\in L^\infty ({\Bbb R})$ such that  
      \begin{equation}
 h= (2\pi)^{-1/2} \Phi\varphi.
\label{eq:pww1}\end{equation}
In this case
        \begin{equation}
       H= \Phi {\bf H} \Phi^* 
      \label{eq:NehC}\end{equation}
      where ${\bf H}$ is the Hankel operator in the space ${\Bbb H}_{+}^2 ({\Bbb R})$ with symbol $\varphi $.
       \end{proposition}

    \begin{pf}   
    Let \e{eq:pww1} hold true.
    Passing to the Fourier transforms, we see that
      \begin{equation}
 {\pmb\la }h, \bar{f} \star f  {\pmb\ra } =\sqrt{2\pi} { \la }\Phi^* h, \bar{w}_\ast
\, w { \ra },  \q   \forall f\in C_{0}^\infty ({\Bbb R}_{+}), \q w=  \Phi^\ast f , \q \forall h\in {\cal S}', 
      \label{eq:Neh2}\end{equation}
      where $w_\ast (\lambda)= w  (-\lambda)$. Under assumption   \e{eq:pww1} the right-hand sides here and in \e{eq:QFR} coincide, and     estimate  \e{eq:QFH1} is true if $\varphi\in L^\infty ({\Bbb R})$. Equality \e{eq:Neh2} implies \e{eq:NehC}.
      
      Conversely, let \e{eq:QFH1} be satisfied so that there exists a bounded operator $H$ satisfying relation \e{eq:QFH}. Define the shift $T(\tau)$, $\tau\geq 0$,  in the space $L^2 ({\Bbb R}_{+})$: $ (T(\tau) f)(t)= f(t-\tau)$. Observe that $ (T(\tau) \bar{f})\star f =\bar{f} \star (T(\tau) f)$ whence $H T(\tau)=T(\tau)^* H$ for all $\tau\geq 0$. It follows that
      $H \Sigma=\Sigma^* H$ where
     \begin{equation}
 \Sigma=-2 \int_{0}^\infty  T(\tau) e^{-\tau} d\tau \q {\rm or} \q    (\Sigma f)(t)= -2 e^{-t}\int_{0}^t e^s f(s) ds.
     \label{eq:Neh2p}\end{equation}
      Passing to the Fourier transforms, we see that the operator ${\bf H}=\Phi^* H \Phi$
     satisfies the commutation relation 
      \begin{equation}
      {\bf H} {\pmb\Sigma}={\pmb\Sigma}^* {\bf H}
        \label{eq:Neh21}\end{equation}
    where according to the second equality \e{eq:Neh2p},  ${\pmb\Sigma}$ is the operator of multiplication by the function $(\lambda-i) (\lambda+i)^{-1}$ in the space ${\Bbb H}_{+}^2 ({\Bbb R})$. Translating the Nehari theorem from the space ${\Bbb H}_{+}^2 ({\Bbb T})$ into ${\Bbb H}_{+}^2 ({\Bbb R})$, we see that  relation  \e{eq:Neh21} implies that ${\bf H}$ is a Hankel operator with symbol $\varphi\in L^\infty ({\Bbb R})$.
       \end{pf}

Finally, we note that the representations in the spaces $l^2 ({\Bbb Z}_{+})$ and $L^2({\Bbb R}_{+})$ can be directly connected by  the unitary operator ${\cal L} $ constructed in terms of the Laguerre functions.  But we not need this construction  here.

 The relations between different representations can be summarized by the following diagrams:  
   \begin{equation}
    \begin{CD}
   u(\mu)   @>>> w(\lambda)= ({\cal U} u )(\nu)
   \\
    @VV    V     @   VV   V     
   \\
    \xi_{n} = ({\cal F} u)_{n}  @>>>  f (t)=(\Phi w) (t) 
    \end{CD}
        \quad\q \q
        \begin{CD}
   {\Bbb H}^2_{+} ({\Bbb T})    @> {\cal U}>>   {\Bbb H}^2_{+} ({\Bbb R})
   \\
    @VV {\cal F}  V     @   VV \Phi V     
   \\
    l^2  ({\Bbb Z_{+}})   @> {\cal L}>>  L^2 ({\Bbb R_{+}})
    \end{CD}  
        \label{eq:dia}\end{equation}
   and
      \begin{equation}
    \begin{CD}
 {\bf  G }   @>>> {\bf  H } = {\cal U}  {\bf  G } {\cal U}^*
   \\
    @VV    V     @   VV   V     
   \\
  G  ={\cal F}  {\bf  G }  {\cal F}^* @>>>  H  =\Phi  {\bf  H }\Phi^*
    \end{CD}
   \q\q \q
   \begin{CD}
   \omega (\mu)   @>>> \varphi(\lambda) 
   \\
    @VV    V     @   VV   V     
   \\
   \varkappa_{n}    @>>>   h (t)      
    \end{CD}
      \label{eq:DIA}\end{equation}
    Of course, the unitary transformations ${\cal F}, {\cal U}, \Phi$ and ${\cal L}$ realizing isomorphisms in \e{eq:dia} are not unique. We can   compose each of them with an automorphism of the corresponding set of Hankel operators  ${\Bbb A}({\cal H})$. 
    The group ${\Bbb G} ({\cal H})$ of automorphisms   of the   set  ${\Bbb A} ({\cal H})$     will be described in the Appendix.

%%%%%%%%%%%%%%%%%%%%%%%%%%%%%%%%%%%%%
 \section{Finite rank Hankel operators   in the   spaces  \\ ${\Bbb H}_{+}^2({\Bbb R})$,       ${\Bbb H}_{+}^2({\Bbb T})$  and  
 $l^2({\Bbb Z}_{+})$}
 %%%%%%%%%%%%%%%%%%%%%%%%%%%%%%%%%%%%%

 Here we reformulate Theorems~\ref{FDH1} and \ref{FDHC} in terms of  Hankel operators acting in the Hardy spaces      ${\Bbb H}_{+}^2({\Bbb R})$ and      ${\Bbb H}_{+}^2({\Bbb T})$ of  analytic functions and in the space of sequences
 $l^2({\Bbb Z}_{+})$. We proceed from relations between various representations   described in Section~4. Now we have to specify  diagrams     \e{eq:DIA} for finite rank Hankel operators.
 
 \medskip  

 {\bf 5.1.}
 Originally, the Kronecker theorem was formulated in the space $l^2({\Bbb Z}_{+})$ (see the paper \cite{Kro} or the book \cite{Pe}, Theorem~3.1 in Chapter~1). It states that,  for a Hankel operator $G$ of rank $r$ determined by the sequence $\varkappa_{n}$, the function
    \begin{equation}
 \omega(\zeta)=\sum_{n=0}^\infty \varkappa_{n} \zeta^n
 \label{eq:pqw}\end{equation}
 is  rational, i.e., $ \omega(\zeta)={\cal P} (\zeta){\cal Q} (\zeta)^{-1}$ where ${\cal P} (\zeta)$ and ${\cal Q} (\zeta)$ are polynomials of degrees $\deg {\cal P}\leq r-1$ and $\deg {\cal Q}\leq r $. Since $\omega\in {\Bbb H}^2_{+}({\Bbb T})$, its poles    lie outside of the unit disc. It follows that for some numbers $\gamma_{m}\in{\Bbb C}$ and $K_{m}\in{\Bbb Z}_{+}$, $m=2,\ldots, M$, function \e{eq:pqw} admits the representation
   \begin{equation}
\omega(\zeta) = R_1 ( \zeta ) +\sum_{m=2}^{M} R_m ( \zeta )   (\zeta- \gamma_m  )^{-K_m-1}, \; |\gamma_{m}|>1, \q R_m (\gamma_m )\neq 0 \; {\rm for} \; m\geq 2,
\label{eq:CH2T}\end{equation}
where all $R_m ( \zeta )$ are polynomials and  $\deg R_m \leq K_m$ for $m\geq 2$. Note that
\[
\rank G = \deg R_{1}+ \sum_{m=2}^M K_m + M.
\]

If  ${\bf G }$ is a finite rank Hankel operator     in the space ${\Bbb H}_{+}^2 ({\Bbb T} )$, then we can apply the Kronecker theorem to the operator $G  ={\cal F} {\bf G } {\cal F}^*$ acting in $l^2({\Bbb Z}_{+})$. Hence    a symbol   $\omega(\zeta)$ of   ${\bf G } $  can be chosen in the form \e{eq:CH2T}.

Similarly, a symbol   $\varphi(z)$ of a finite rank Hankel operator  ${\bf H } $  in the space ${\Bbb H}_{+}^2 ({\Bbb R} )$ can be chosen in the form
 \begin{equation}
\varphi(z) =\sum_{m=1}^M  Q_m (z) ( \alpha_m -i z)^{-K_m-1}, \q \Re \alpha_m>0, \q Q_m (-i\alpha_m)\neq 0,
\label{eq:CH2x}\end{equation}
where $Q_m (z)$ are polynomials and  $\deg Q_m \leq K_m$. The relation between symbols \e{eq:CH2T} and \e{eq:CH2x}  of the operators ${\bf G} $  and $   {\bf H} = {\cal U}  {\bf G} {\cal U}^*  $ is given by formula \e{eq:pw}. In particular, we have
\[
\alpha_{m}= \frac{ \gamma_{m}+1}{  \gamma_{m} -1}
\] 
and
 \begin{equation}
  \begin{split}
 Q_1 (-i )& = (-1)^{K_1} 2^{K_{1}+1}\frac{1}{K_{1}!}R_1^{(K_{1})} ,\q  K_1 =\deg R_1,   \q \alpha_{1}= 1  ,
 \\
 Q_m (-i\alpha_m) &= -2^{K_{m}+1} \gamma_{m} (\gamma_{m}-1)^{-2K_{m}-2}R_m (\gamma_m ),\q \alpha_{m}\neq 1.
  \end{split}
\label{eq:CL}\end{equation}

Since
\[
k!\int_{-\infty}^\infty (\alpha-iz)^{-k-1} e^{-izt} d z=2\pi t^k e^{-\alpha t},\q \Re\alpha>0,
\]
the operator $H=\Phi {\bf H}\Phi^*$ acts by formula  \e{eq:H1} and has kernel \e{eq:FDvm} where $P_{m}(t)$ is a polynomial of degree $K_{m}$ and
 \begin{equation}
 P_{m}^{(K_{m})} = Q_{m} (-i\alpha_{m}).
\label{eq:CL1}\end{equation}
Thus, as stated in the Introduction, all finite rank Hankel operators in the space $L^2 ({\Bbb R}_{+})$ are given by 
formulas  \e{eq:H1}, \e{eq:FDvm}. 

Finally,     for a finite rank Hankel operator $G $ in the   space $l^2 
   ({\Bbb Z}_{+})$, the operator ${\bf G}={\cal F}^* G {\cal F} $ has symbol   \e{eq:CH2T}. Expanding this function   into the Fourier series and using relation  \e{eq:pww}, we find     that
  \begin{equation}
\varkappa_{n}= \tau_{n}  + \sum_{m=2}^M   T_m ( n) q_{m}^n , \q |q_m| <1, 
\label{eq:CH2Z}\end{equation}
where $\tau_{n}  = 0$  for   $n\geq K_1+1$ 
and $T_{m}$ are polynomials of degree $K_{m}$. We note that $q_{m}=\gamma_{m}^{-1}$ and
   \begin{equation}
 K_{1} !  \tau_{K_{1}}= R_{1}^{(K_{1})} , \q T_{m}^{(K_{m})} =   (-1)^{K_{m}+1} q_{m}^{K_{m}+1}  R(\gamma_{m}).
\label{eq:CL2}\end{equation}

Let us now consider the self-adjoint case.  If  ${\bf G}={\bf G}^*$   and the sum in \e{eq:CH2T} contains a term with $ \gamma_m $, then necessarily it also contains the term  with   $\bar{  \gamma}_m $. We  suppose that $\Im \gamma_{m}=0$ for $m=2,\ldots, M_{0}$ and $\Im \gamma_{m}< 0$,  $\gamma_{ M_{1}+m}=\bar{\gamma}_{ m}$ for $m=M_{0} +1, \ldots,M_{0} + M_{1} $. 
Then  $\ov{R_{m} (\bar{\zeta })}=R_{m} ( \zeta)$ for $m=1,\ldots, M_{0}$ and  $\ov{R_{M_1+m} (\bar{\zeta} )}=R_{m} ( \zeta )$   for $m=M_{0} +1, \ldots,M_{0} + M_{1} $. 

Similarly, if  ${\bf H} ={\bf H}^*$  and the sum in \e{eq:CH2x} contains a term with $ \alpha_m $, then necessarily it also contains the term with   $\bar{ \alpha}_m $. We again suppose that $\Im \alpha_{m}=0$ for $m=1,\ldots, M_{0}$ and $\Im \alpha_{m}>0$,  $\alpha_{ M_{1}+m}=\bar{\alpha}_{ m}$ for $m=M_{0} +1, \ldots,M_{0} + M_{1} $.  Then  $\ov{Q_{m} (\bar{z} )}=Q_{m} ( -z )$ for $m=1,\ldots, M_{0}$ and  $\ov{Q_{M_1+m} (\bar{z} )}=Q_{m} ( -z )$   for $m=M_{0} +1, \ldots,M_{0} + M_{1} $. 

 Finally, if $G =G^*$, then necessarily $\tau_{n}  =\bar{\tau}_{n}  $ and  if the sum in \e{eq:CH2Z} contains a term $T_m ( n) q_{m}^n$, then it also contains the term $\ov{T_m ( n)}\bar{q}_{m}^n$.  We  suppose that $\Im q_{m}=0$ for $m=2,\ldots, M_{0}$, $\Im q_{m}>0$ for $m=M_{0}+ 1,\ldots, M_{0}+M_{1}$ and $q_{ M_{1}+m}=\bar{q}_{ m}$ for $m=M_{0} +1, \ldots,M_{0} + M_{1} $. 
The coefficients of the polynomials  $T_{m} $ for $m=2,\ldots, M_{0}$ are of course real.

\medskip  

 {\bf 5.2.}
  Now we are in a position to reformulate Theorems~\ref{FDH1} and \ref{FDHC} in various representations of Hankel operators. Recall that the numbers $ {\sf p}_{m}$ were defined by formula  \e{eq:ppq}.  Let us start with finite rank Hankel operators $\bf H$ in the space ${\Bbb H}_{+}^2({\Bbb R})$. 
  Given the one-to-one correspondence between Hankel operators with kernels  \e{eq:FDvm} and symbols \e{eq:CH2x} and, in particular, equality \e{eq:CL1}, the following result is equivalent to Theorem~\ref{FDH1}.
   
  \begin{theorem}\label{FDHR}
 Let  the symbol of a self-adjoint Hankel operator  $\bf H$  in the space ${\Bbb H}_{+}^2 ({\Bbb R}) $ be given by formula \e{eq:CH2x} where  $Q_m (z)$ are polynomials of degree $\deg Q_m \leq K_m$. Let the numbers ${\cal N}_\pm^{(m)}$ be defined by formula \e{eq:RXm} where ${\sf p}_{m}= Q_{m} ( -i \alpha_m )$.    Then the total numbers $N_{\pm} ({\bf H})$  of   $($strictly$)$ positive and negative eigenvalues of the operator $ {\bf H}$ are given by the formula 
   \e{eq:TN}.
 \end{theorem}

     In particular (cf. Corollary~\ref{FD}),  ${\bf H }\geq 0$ (${\bf H }\leq 0$) if and only if all poles of its symbol lie on the imaginary axis, are simple, the real parts of the residues are equal to zero and their imaginary parts are positive (negative).

    Note also that according to \e{eq:pww1} the symbol of the Carleman operator can be chosen as $\varphi_{0}(\lambda)=\pi i\sgn \lambda$. Therefore the next result is a direct consequence of Theorem~\ref{FDHC}.

  \begin{theorem}\label{FDHRk}
    Let $\bf H$ be the Hankel operator with symbol  $  \pi i\sgn \lambda+\varphi (\lambda)$ where $\varphi (\lambda)$ is function \e{eq:CH2x}.
  Then  the total number $N_{-} ({\bf H})$  of     its negative eigenvalues    is   given by formula   \e{eq:TN}.
 \end{theorem}

Quite similarly, 
  given the one-to-one correspondence between Hankel operators with   symbols \e{eq:CH2T} and  \e{eq:CH2x} and, in particular, equalities \e{eq:CL}, the following result is equivalent to Theorem~\ref{FDHR}.

  \begin{theorem}\label{FDHT}
 Let  the symbol of a self-adjoint Hankel operator  $\bf G$  in the space ${\Bbb H}_{+}^2 ({\Bbb T}) $ be given by formula \e{eq:CH2T}  where  $R_m (\zeta)$ are polynomials of degree $\deg R_m \leq K_m$.
 Let the numbers ${\cal N}_\pm^{(m)}$ be defined by formula \e{eq:RXm} where
 \[  
 {\sf p}_1=  R_1^{(K_1)}   \q {\rm and} \q  {\sf p}_m= -  R_m (\gamma_{m}) \sgn \gamma_{m}  \q {\rm if} \q m=2,\ldots, M_0.
 \]
   Then the total numbers $N_{\pm} ({\bf G})$  of   $($strictly$)$ positive  and negative eigenvalues of the operator $ {\bf G}$ are given by the formula    \e{eq:TN} $($if $R_{1}(\zeta)=0$, then the first sum in \e{eq:TN}  starts with $m=2)$.
 \end{theorem}   
   
     In particular (cf.   Corollary~\ref{FD}),  ${\bf G }\geq 0$  ( ${\bf G }\leq 0$) if and only if all poles of its symbol lie on the real axis, are simple and the residues are positive (negative); moreover, it is required that $\deg R_{1}=0$ and $R_1 \geq 0$ ($R_1 \leq 0$).
     
     Theorem~\ref{FDHRk} can also be  reformulated  in an obvious  way in terms of Hankel operators in the space ${\Bbb H}_{+}^2({\Bbb T})$.

  \begin{theorem}\label{FDHRk1}
    Let $\bf G$ be the Hankel operator with symbol  $  \pi i \mu^{-1} \sgn \Im \mu+ \omega(\mu)$ where $\omega (\mu)$ is function \e{eq:CH2T}.
  Then  the total number $N_{-} ({\bf G})$  of     its negative eigenvalues    is   given by formula   \e{eq:TN}.
 \end{theorem}  

Finally, we use
  the one-to-one correspondence between Hankel operators with   symbols   \e{eq:CH2T} 
  and with matrix elements \e{eq:CH2Z} and, in particular, equalities \e{eq:CL2}. Therefore the following result is equivalent to Theorem~\ref{FDHT}.

 \begin{theorem}\label{FDHSe}
 Let    $G$ be a finite rank Hankel operator in the space
    $L^2 ({\Bbb Z}_+) $ with matrix elements \e{eq:CH2Z} where $ \tau_n =0$ for $n> K_{1}$, $ \tau_{K_1 }\neq 0$ and  $T_m  $ are polynomials of degree $  K_m$.  Let the numbers ${\cal N}_\pm^{(m)}$ be defined by formula \e{eq:RXm} where  
 \[  
  {\sf p}_1=  \tau_{K_1 } \q {\rm and} \q
   {\sf p}_m=   T_m^{(K_m )} \q {\rm if }\q m=2,\ldots, M_0.  
 \]
    Then the total numbers $N_{\pm} (G)$  of   $($strictly$)$ positive  and negative eigenvalues of the operator $ G$ are given by the formula    \e{eq:TN} $($if $\tau_{n}=0$ for all $n \geq 0$, then the first sum in \e{eq:TN}  starts with $m=2)$.
 \end{theorem}

     In particular (cf.   Corollary~\ref{FD}),  $G \geq 0$  ( $G \leq 0$) if and only if 
      \[
\varkappa_{n}= t_1\d_{n,0}+ \sum_{m=2}^{M_0}   t_m   q_{m}^n , \q q_m \in (-1,1),
\]
where all numbers $t_1, \ldots, t_{M_0}$ are positive (negative).

     Theorem~\ref{FDHRk1} can also be  reformulated  in an obvious  way in terms of operators in the space   $l^2({\Bbb Z}_+) $ if one takes into account that the matrix elements of the Carleman operator equal $\varkappa_n^{(0)}=2(n+1)^{-1}$ for $n$ even and $\varkappa_n^{(0)} = 0$ for $n$ odd.
     
      \begin{theorem}\label{FDHRk2}
    Let $  G$ be the Hankel operator  in the space
   $l^2({\Bbb Z}_+) $ with matrix elements $\varkappa_n^{(0)} + \varkappa_{n}$ where the numbers $\varkappa_{n}$  are defined by formula \e{eq:CH2Z}. 
  Then  the total number $N_{-} ( G )$  of     its negative eigenvalues    is   given by formula   \e{eq:TN}.
 \end{theorem}

 \appendix{}

\section{The automorphism group of Hankel operators}

 {\bf A.1.}
 Let $\cal H$ be one of the spaces ${\Bbb H}_{+}^2 ({\Bbb R})$,  $L^2 ({\Bbb R}_{+})$, 
  ${\Bbb H}_{+}^2 ({\Bbb T})$ or $l^2 ({\Bbb Z}_{+})$. 
 Our goal here is to describe the group ${\Bbb G} ({\cal H})$ of all automorphisms of the set $ {\Bbb A}({\cal H})$ of Hankel operators in ${\cal H}$.
        By definition, a unitary operator $U\in {\Bbb G} ({\cal H})$ if and only if $U H U^*\in {\Bbb A}({\cal H})$ for all $H\in {\Bbb A}({\cal H})$.
      Of course, for a Hankel operator $H$ and an arbitrary unitary operator $U$, the operator $U H U^*$ is not necessarily Hankel. Hence the group ${\Bbb G} ({\cal H})$ is smaller than the group of all unitary operators. It turns out that this group admits a simple description.

    It is sufficient to describe ${\Bbb G} ({\cal H})$ for one of the spaces 
${\Bbb H}_{+}^2 ({\Bbb R})$, $L^2 ({\Bbb R}_{+})$,  ${\Bbb H}_{+}^2 ({\Bbb T})$ or $  l^2 ( {\Bbb Z}_{+})$. We choose ${\cal H} = {\Bbb H}_{+}^2 ({\Bbb R})$.
Then other groups are obtained by conjugations with  the unitary transformations $\Phi, {\cal U}^*$ and $ {\cal F}$:
 \begin{align*}
{\Bbb G} (L^2 ({\Bbb R}_{+}))&=\Phi {\Bbb G} ({\Bbb H}_{+}^2 ({\Bbb R}))\Phi^*,
\\
{\Bbb G} ({\Bbb H}_{+}^2 ({\Bbb T}) )&={\cal U}^* {\Bbb G} ({\Bbb H}_{+}^2 ({\Bbb R})){\cal U},
\\
{\Bbb G} (l^2 ( {\Bbb Z}_{+}) )&={\cal F}  {\Bbb G} ({\Bbb H}_{+}^2 ({\Bbb T})) {\cal F}^*.
 \end{align*}

 Let us define the dilation operators ${\bf D}_\rho$, $\rho>0$,  in the space  $ {\Bbb H}_{+}^2 ({\Bbb R})$:
 \[
( {\bf D}_\rho u)(\lambda) = \rho^{1/2} u(\rho \lambda).
\]
Obviously,   the operators 
  ${\bf D}_\rho$    are unitary.  
Set
 \[
(\pmb{\cal I} u)(\lambda) =i \lambda^{-1} u(-\lambda^{-1}).
\]
Then $\pmb{\cal I}  : {\Bbb H}_{+}^2 ({\Bbb R})\to {\Bbb H}_{+}^2 ({\Bbb R})$, $\pmb{\cal I} $ is the involution, i.e. $\pmb{\cal I}  ={\pmb{\cal I}}^2$, and $\pmb{\cal I} $ is also unitary. It is  easy to see that 
      \begin{equation}
 {\bf D}_{\rho}{\bf H} (\varphi)  {\bf D}_{\rho}^*={\bf H} (\varphi_{\rho}) \q  {\rm and}\q
\pmb{\cal I}  {\bf H} (\varphi)  {\pmb{\cal I} }^* ={\bf H} (\ti{\varphi}  )
\label{eq:Aut9}\end{equation}
where
$\varphi_{\rho}(\lambda)=\varphi   (\rho\lambda)$ and $\ti{\varphi}(\lambda)=\varphi(-\lambda^{-1})$. In particular,
$  {\bf D}_\rho \in {\Bbb G} ({\Bbb H}_{+}^2 ({\Bbb R}))$ and $\pmb{\cal I}  \in {\Bbb G} ({\Bbb H}_{+}^2 ({\Bbb R}))$. It turns out that the group $ {\Bbb G} ({\Bbb H}_{+}^2 ({\Bbb R}))$ is exhausted by these transformations.
Let us state the precise result.

 \begin{theorem}\label{Aut}
A unitary operator ${\bf U}\in  {\Bbb G} ({\Bbb H}_{+}^2 ({\Bbb R}))$ if and only if it has one of the two forms:  ${\bf U}=\theta {\bf D}_\rho$ or ${\bf U}=\theta {\bf D}_\rho \pmb{\cal I} $ for some $\theta\in {\Bbb T}$ and $\rho>0$.
   \end{theorem} 
   
   Actually, we shall prove a stronger statement.
   
    \begin{theorem}\label{Aut1}
    Let ${\bf H}_{\alpha}$ be the Hankel operator in the space ${\Bbb H}_{+}^2 ({\Bbb R})$ with symbol $\varphi_{\alpha}(\lambda)= 2\alpha (\alpha-i\lambda)^{-1}$. Suppose that an operator ${\bf U}$ is  unitary and ${\bf U} {\bf H}_{\alpha} {\bf U}^* \in {\Bbb A} ({\Bbb H}_{+}^2 ({\Bbb R}))$ for all $\alpha>0$. Then either ${\bf U} =\theta {\bf D}_\rho$ or ${\bf U} =\theta {\bf D}_\rho\pmb{\cal I}  $ for some $\theta\in {\Bbb T}$ and $\rho>0$.
   \end{theorem} 

 \begin{pf}
 Set $U=\Phi {\bf U} \Phi^*$ and  $H_{\alpha}=\Phi {\bf H}_{\alpha} \Phi^* $. It follows from formula \e{eq:pww1} that $H_\alpha$ is the Hankel operator in the space $L^2 ({\Bbb R}_{+})$ with kernel $h_\alpha (t)=2 \alpha e^{- \alpha t}$, that is, $H_\alpha f= (f,\psi_\alpha) \psi_\alpha$ where
     \begin{equation}
     \psi_\alpha(t) =\sqrt{2 \alpha}  e^{- \alpha t}.
\label{eq:Aut2} \end{equation}
   By our assumption, the operator 
 $U H_\alpha U^* \in {\Bbb A}(L^2 ({\Bbb R}_{+}))$.  It has rank one, and its non-zero eigenvalue equals $1$. By the Kronecker theorem, all rank one
   Hankel operators have kernels $p e^{-\betaÊt}$ for some   $p, \beta\in {\Bbb C}$ with $\Re\beta>0$. They are self-adjoint and have the eigenvalue   $1$ if and only if 
      $\beta>0$ and  $p=\sqrt{2\beta}$.  Therefore 
 \[
  U H_\alpha U^*= H_\beta 
  \]
  and hence
  \[
  (f, U \psi_\alpha ) U\psi_\alpha =  (f, \psi_\beta )\psi_\beta 
  \]
  for all $f\in L^2 ({\Bbb R}_{+})$ and some $\beta= v(\alpha)$. It follows that 
   \begin{equation}
 U\psi_\alpha=\theta (\alpha) \psi_{v(\alpha)},\q |\theta (\alpha)|=1.
\label{eq:Aut}
\end{equation}

We have to find the functions   $\theta (\alpha)$ and $v(\alpha)$.
  Let us   take the unitarity of $U$ into account.  Since   $(U \psi_{\alpha_{1}}, U \psi_{\alpha_{2}})=      ( \psi_{\alpha_{1}},   \psi_{\alpha_{2}})$,  relation \e{eq:Aut} implies that 
       \begin{equation}
\theta (\alpha_{1})\ov{\theta (\alpha_2)}( \psi_{v(\alpha_{1})},   \psi_{v(\alpha_{2})})=( \psi_{\alpha_{1}},   \psi_{\alpha_{2}}), \q \forall \alpha_{1}, \alpha_{2} >0.
\label{eq:Aut1}\end{equation}
Note that $\psi_{\alpha}(t)>0$ for all $\alpha>0$ and $t>0$ and hence $\theta (\alpha_{1})\ov{\theta (\alpha_2)}>0$. Using also that $|\theta (\alpha)|=1$, we see that $\theta (\alpha )= \theta (1)$ for all $\alpha >0$; thus $\theta (\alpha)=: \theta$ does not depend on $\alpha$. Returning to \e{eq:Aut1} and using the explicit expression \e{eq:Aut2} for $\psi_{\alpha} (t)$, we obtain the equation
  \begin{equation}
  \frac{\sqrt{v(\alpha_{1})v(\alpha_{2})}}{v(\alpha_{1})+v(\alpha_{2})}
  =  \frac{\sqrt{ \alpha_{1} \alpha_{2} }}{\alpha_{1} +\alpha_{2} }.
\label{eq:Aut3}
\end{equation}

Set here $\alpha_{1}=1$,  $\alpha_{2}= \alpha$. Equation \e{eq:Aut3} for $v(\alpha)$ has two solutions 
  \[
 v( \alpha )=    v(1 )\alpha \q {\rm and}  \q  v(\alpha )=    v(1 )\alpha^{-1}.
\]
It now follows from  \e{eq:Aut} that
 \[
 U\psi_\alpha=\theta   \psi_{\rho^{-1} \alpha} \q {\rm or}  \q
  U\psi_\alpha=\theta   \psi_{(\rho \alpha)^{-1}}
\]
where $\rho=v(1)^{-1}$. Since $\Phi\varphi_{\alpha}=2 \sqrt{\pi\alpha}\psi_{\alpha}$, these equalities can be rewritten as
 \begin{equation}
 {\bf U}\varphi_\alpha=\theta \sqrt{\rho} \varphi_{\rho^{-1} \alpha} \q {\rm or}  \q
  {\bf U}\varphi_\alpha =\theta \alpha \sqrt{\rho} \varphi_{(\rho \alpha)^{-1}} .
\label{eq:Aut5y}\end{equation}
Note that
 \[
 {\bf D}_{\rho}\varphi_\alpha=  \sqrt{\rho} \varphi_{\rho^{-1} \alpha} \q {\rm and}  \q
\pmb{\cal I} \varphi_\alpha =  \alpha   \varphi_{  \alpha^{-1}}.
\]
Hence  \e{eq:Aut5y} are equivalent to the equalities 
 \[
 {\bf U}\varphi_\alpha= \theta  {\bf D}_{\rho}\varphi_\alpha \q {\rm or}  \q
 {\bf U}\varphi_\alpha=  \theta  {\bf D}_{\rho} \pmb{\cal I} \varphi_{  \alpha},\q \forall \alpha>0 .
\]

It remains to extend these relations to the whole space ${\Bbb H}_{+}^2 ({\Bbb R})$. To that end, we have to show that the set of the functions $\varphi_{\alpha}  $ where $\alpha>0$ is arbitrary is dense in ${\Bbb H}_{+}^2 ({\Bbb R})$ or, equivalently, that
the set of the functions $\psi_{\alpha}  $    is dense in $L^2 ({\Bbb R}_{+})$.
Set 
\[
(Lf)(t)=\int_{0}^\infty e^{-\alpha t} f(\alpha) d\alpha.
\]
This operator is self-adjoint and bounded in the space  $L^2 ({\Bbb R}_{+})$. It has purely absolutely continuous spectrum (see, e.g., \cite{Y3}). Therefore its range is dense in $L^2 ({\Bbb R}_{+})$.
   \end{pf}

    \medskip  

 {\bf A.2.}
 Let us now describe the group ${\Bbb G} ({\cal H})$ in other representations of the space ${\cal H}$. In view of Theorem~\ref{Aut}, to that end we only have to calculate the operators 
    \[
 {\bf W}_{\rho}= {\cal U}^*  {\bf D}_{\rho} \,{\cal U}, \q  {D}_{\rho}= \Phi  {\bf D}_{\rho}\Phi^*,\q
 {W}_{\rho}= {\cal F}  {\bf W}_{\rho}{\cal F}^*
\]
 and
   \[
\pmb{\cal J} = {\cal U}^* \pmb{\cal I} {\cal U}, \q  {\cal I}= \Phi \pmb{\cal I}\Phi^*,\q
{\cal J} = {\cal F} \pmb{\cal J} {\cal F}^*
\]
acting in the spaces ${\Bbb H}^2({\Bbb T})$, $L^2({\Bbb R}_{+})$, $l^2 ({\Bbb Z}_{+})$, respectively. 

According to formula \e{eq:pwtq}   we have
\[
( {\bf W}_{\rho} f)(\mu)= \tfrac{2 \rho^{1/2}}{ \rho+1} \tfrac{1}{  \tau(\rho)\mu+1}
 f(\tfrac{\mu+\tau(\rho)}{ \tau(\rho)\mu+1}), \q\tau(\rho)=
\tfrac{ \rho-1}{  \rho +1}\in (-1,1),
\]
and
 \begin{equation}
( \pmb{\cal J} f)(\mu)=f(-\mu).
\label{eq:Aut7}\end{equation}
The role of \e{eq:Aut9} is now played by the relations
   \[
  {\bf W}_{\rho}{\bf G} (\omega)   {\bf W}_{\rho}^* ={\bf G} (\omega_{\rho}) \q  {\rm and}\q
 \pmb{\cal J} {\bf G} (\omega)   \pmb{\cal J}^*={\bf G} (\ti{\omega} )
\]
where
\[
\omega_{\rho}(\mu)= \omega ( \tfrac{\mu+\tau(\rho)}{ \tau(\rho)\mu+1}) \q  {\rm and}\q
\ti{\omega}(\mu)=\omega(-\mu).
\]

The operator $ {D}_{\rho}$ is again the dilation, $( {D}_{\rho} f)(t)= \rho^{-1/2} f(\rho^{-1}t)$, and for a Hankel operator $H$ with kernel $h(t)$, the operator
$ {D}_{\rho} H  {D}_{\rho}^*$  has kernel $h_{\rho}(t)= h(\rho^{-1} t)$. Apparently there is no simple formula for the operator $ {\cal I}$.

It follows from \e{eq:Aut7} that  $({\cal J} \xi)_{n} = (-1)^n \xi_{n}$, and 
  for a Hankel operator $G$ with matrix elements $\varkappa_{n}$, the operator
${\cal J} G {\cal J}^*$  has matrix elements $(-1)^n \varkappa_{n}$. On the contrary, there seems to be no direct expression for the operators $ {W}_{\rho}$.

%%%%%%%%%%%%%%%%%%%%%%%%%%%%%%%%%%%%%%%%
%%%%%%%%%%%%%%%%%%%%%%%%%%%%%%%%%%%%%%%%

 \end{document}